\documentclass{amsart}

\usepackage{graphicx}
\usepackage{a4wide}
\usepackage{psfrag}
\usepackage[all]{xy}

\let\cal\mathcal
\def\AA{{\cal A}}
\def\BB{{\cal B}}
\def\CC{{\cal C}}

\def\FF{{\cal F}}

\def\KK{{\cal K}}
\def\LL{{\cal L}}

\def\PP{{\cal P}}
\def\QQ{{\cal Q}}

\def\TT{{\cal T}}
\def\UU{{\cal U}}
\def\VV{{\cal V}}
\def\WW{{\cal W}}

\let\blb\mathbb

\def\bZ{{\blb Z}}

\def\bN{{\blb N}}

\def\bZ{{\blb Z}}
\def\bZz{\bZ_{0}}

\let\frak\mathfrak

\def\Mod{\operatorname{Mod}}

\def\mod{\operatorname{mod}}

\def\rad{\operatorname {rad}}

\def\Rep{\operatorname {Rep}}

\def\repc{\rep^{\rm{cfp}}}
\def\rep{\operatorname{rep}}

\def\Ext{\operatorname {Ext}}
\def\Hom{\operatorname {Hom}}

\def\End{\operatorname {End}}

\def\cone{\operatorname {cone}}
\def\coker{\operatorname {coker}}

\def\End{\operatorname {End}}

\def\add{\operatorname {ind}}

\def\infrad{\operatorname {rad}^{\infty}}
\def\radn{\operatorname {rad}^{n}}

\def\r{\rightarrow}

\def\card{\operatorname {card}}

\DeclareMathOperator{\Ind}{Ind}

\DeclareMathOperator{\Aut}{Aut}

\DeclareMathOperator{\ind}{ind}

\newcommand\phic{\phi^\text{comp}}
\newcommand\phio{\phi^\text{obj}}

\newcommand\TXZ{\TT \stackrel{\rightarrow}{\times} \bZ}
\newcommand\Db{D^{b}}
\newcommand\Kb{K^{b}}
\newcommand\AL{A_{\LL}}
\newcommand\DL{D_{\LL}}

\newcommand\tri[3]{#1\to #2\to #3\to #1[1]}
\newcommand\ZAi{\bZ A_{\infty}}
\newcommand\ZAii{\bZ A_{\infty}^{\infty}}
\newcommand\ZDi{\bZ D_{\infty}}
\newcommand\asa{\Leftrightarrow}
\newcommand\St{S_{\to}}
\newcommand\Tt{T_{\to}}

\newtheorem{lemma}{Lemma}[section]
\newtheorem{proposition}[lemma]{Proposition}
\newtheorem{theorem}[lemma]{Theorem}

\theoremstyle{definition}

\newtheorem{example}[lemma]{Example}
\newtheorem{definition}[lemma]{Definition}

{

}

\theoremstyle{remark}

\newtheorem{remark}[lemma]{Remark}

\newdimen\uboxsep \uboxsep=1ex
\def\uboxn#1{\vtop to 0pt{\hrule height 0pt depth 0pt\vskip\uboxsep
\hbox to 0pt{\hss #1\hss}\vss}}

\def\uboxs#1{\vbox to 0pt{\vss\hbox to 0pt{\hss #1\hss}
\vskip\uboxsep\hrule height 0pt depth 0pt}}

\def\Ob{\operatorname{Ob}}

\newcommand\exa{\nopagebreak \begin{center}\smallskip \nopagebreak               \begin{minipage}[t]{6cm}\sloppy}
\newcommand\exb{\end{minipage}\kern 1cm\begin{minipage}[t]{8cm}\sloppy}
\newcommand\exc{\end{minipage}\kern -3cm \smallskip\end{center}}

\begin{document}

\title{Classification of abelian hereditary directed categories satisfying Serre duality}
\author{Adam-Christiaan van Roosmalen}
\address{Adam-Christiaan van Roosmalen\\Hasselt University
\\Research group Algebra\\Agoralaan, gebouw D\\B-3590 Diepenbeek (Belgium)}\email{AdamChristiaan.vanRoosmalen@UHasselt.be}
\subjclass[2000]{16G20; 16G70; 18E10; 18E30}

\begin{abstract}
In an ongoing project to classify all hereditary abelian categories, we provide a classification of $\Ext$-finite directed hereditary abelian categories satisfying Serre duality up to derived equivalence.

In order to prove the classification, we will study the shapes of the Auslander-Reiten components extensively and use appropriate generalizations of tilting objects and coordinates, namely partial tilting sets and probing of objects by quasi-simples. 
\end{abstract}

\maketitle

\tableofcontents

\section{Introduction}

Let $k$ be an arbitrary algebraically closed field. In this paper we classify, up to
derived equivalence, $k$-linear $\Ext$-finite \emph{directed
  hereditary categories satisfying Serre duality} (see below for
definitions).  In this way we accomplish a step in the ongoing
classification project of hereditary categories. As an $\Ext$-finite
hereditary category may be viewed as a homological generalization of a
smooth projective curve, one of the motivations for this
classification project is non-commutative algebraic geometry (see e.g.\ \cite{Stafford01}).

\medskip
Our classification is a natural complement to the following results.
\begin{enumerate}
\item In \cite{Happel01} Happel classifies $\Ext$-finite hereditary
  categories with a tilting object.
\item\label{enumerate:ReVdB02} In \cite{ReVdB02} Reiten and Van den Bergh classify $\Ext$-finite
  Noetherian hereditary categories satisfying Serre duality.
\end{enumerate}
The hereditary categories considered in (\ref{enumerate:ReVdB02}) fall naturally into several classes. One class contains mild generalizations of smooth projective curves and the other classes are in some way associated to quivers.

\medskip

In an early version of \cite{ReVdB02} it was conjectured that the 
hereditary categories appearing in \cite{ReVdB02} constitute in fact a
complete classification of $\Ext$-finite hereditary categories with
Serre duality \emph{up to derived equivalence}.  This conjecture was
quickly shown to be false by Ringel \cite{Ringel02} who constructed
counterexamples.  Since these counterexamples are directed it was
very natural for us to try to classify directed hereditary categories.

\medskip

To precisely state our classification result we now give some definitions. Let
$\AA$ be a $k$-linear abelian category.  
\begin{enumerate}
\item
We say that $\AA$ is \emph{hereditary} if $\Ext^2_\AA(X,Y)=0$ for $X$, $Y$ arbitrary
objects in $\AA$.
\item
We say that $\AA$ is \emph{$\Ext$-finite} if $\dim \Ext^i_{\AA}(X,Y)<\infty$
for all $i$.
\item
Assume  that $\AA$ is $\Ext$-finite. We say that $\AA$ satisfies \emph{Serre duality} \cite{Bondal89} if there is an auto-equivalence $F$ of $D^b(\AA)$ such that
there are isomorphisms of $k$-vector spaces
\[
\Hom_{D^b(\AA)}(X,Y)\cong \Hom_{D^b(\AA)}(Y,FX)^\ast
\]
natural in $X,Y$. If $\AA$ is hereditary then if $\AA$ has Serre duality it has
almost split sequences, and the converse is true if $\AA$ has neither
 projectives nor injectives \cite{ReVdB02}.
\item Assume that $\AA$ is a Krull-Schmidt category (e.g.\ if $\AA$ is
$\Hom$-finite). Then $\AA$ is \emph{directed} if there is no cycle
of maps between indecomposable objects
\[
X_0\r X_1\r \cdots \r X_n\r X_0
\]
which are neither zero, nor isomorphisms. 

It follows in particular that $\Ext^1_\AA(X,X)=0$ for any indecomposable object $X$. Thus any indecomposable object is rigid.
In fact: experience suggests that the entire structure of 
directed categories is very rigid and that their structure is entirely controlled by
combinatorics.
\end{enumerate}
Our classification will be stated in terms of the representation
theory of certain partially ordered sets. If $\PP$ is a partially
ordered set then we may view it as a category where an arrow between
objects $x,y\in \PP$ exists if and only if $x<y$. The category
$\Rep(\PP)$ of \emph{$\PP$-representations} is the category of covariant functors from $\PP$ to $k$-vector spaces. 

To every element of $\PP$ there  are naturally associated
an indecomposable projective object as well as an indecomposable injective object. We call such objects \emph{standard} projective and injective objects. We say that
a representation is \emph{finitely presented} if it has a finite
presentation by finite direct sums of standard projectives.
\emph{Cofinitely presented} is defined dually. We define $\repc(\PP)$
as the full subcategory of $\Rep(\PP)$ consisting of
representations which are both finitely and cofinitely presented. The following
is a special case of Proposition \ref{proposition:ShIsLocal}.
\begin{proposition} Assume that $\PP$ is a \emph{forest} (i.e.\ $\PP$
does not contain $x,y,z,t$ such that $x<y<t$ and $x<z<t$ but
$y$, $z$ are incomparable). Then $\rep(\PP)$ is a hereditary abelian
category.
\end{proposition}

Below $\LL$ will be a totally ordered ordered set in which
every element has an immediate successor and an immediate predecessor.  It is easy to see
that any such partially ordered set is of the form
$\TT\overset{\r}{\times} \bZ$ where $\TT$ is totally ordered and
$\overset{\r}{\times}$ denotes the lexicographically ordered product.
Pictorially we may draw $\LL$ as
\[
\cdots [\cdots \r\bullet \r \bullet\r\cdots \r\bullet\r \bullet\r\cdots]\cdots [\cdots \r\bullet \r \bullet\r\cdots \r\bullet \r \bullet\r\cdots]\cdots
\]
If $\LL=\bZ$ then $\LL$ is sometimes referred to as a $A_{\infty}^\infty$
quiver. Therefore we will usually write $A_\LL$ for $\LL$. We also
define $D_\LL$ as the union of $A_\LL$ with two distinguished objects
 which are strictly smaller than the elements of $A_\LL$ but
incomparable with each other. Schematically:
\[
\hskip -12cm\xymatrix@=0.2cm{
\bullet\ar[dr]&\\
&\hbox to 0pt{$\cdots [\cdots \r\bullet \r \bullet\r\cdots \r\bullet\r \bullet\r\cdots]\cdots [\cdots \r\bullet \r \bullet\r\cdots \r\bullet \r \bullet\r\cdots]\cdots
$\hss}\\
\bullet\ar[ur]&
}
\]
The following is our main result (Theorem \ref{theorem:MainTheorem} in the text).
\begin{theorem} \label{introtheorem} A connected directed hereditary category $\AA$ satisfying
Serre duality is derived equivalent to $\repc(\PP)$ where $\PP$ is either
a Dynkin quiver, $A_\LL$ or $D_\LL$.
\end{theorem}
The  categories $\repc(A_\LL)$ with $\LL=\TT\overset{\r}{\times} \bZ$
for $|\TT|>2$ were precisely the above mentioned counterexamples constructed by Ringel in \cite{Ringel02}. 

\medskip

The categories occurring in Theorem \ref{introtheorem} have rather attractive Auslander-Reiten quivers (see \S\ref{section:Sh}).  If $\PP=A_\LL$ then the Auslander-Reiten quivers of $\repc(A_\LL)$ and its derived category have the form:

\vspace*{2mm}

\begin{center}
\psfrag{T}[][]{$\LL$}
\psfrag{A}[][]{$\AA$}
\psfrag{a1}[][]{$\AA[1]$}
\includegraphics[height=5cm]{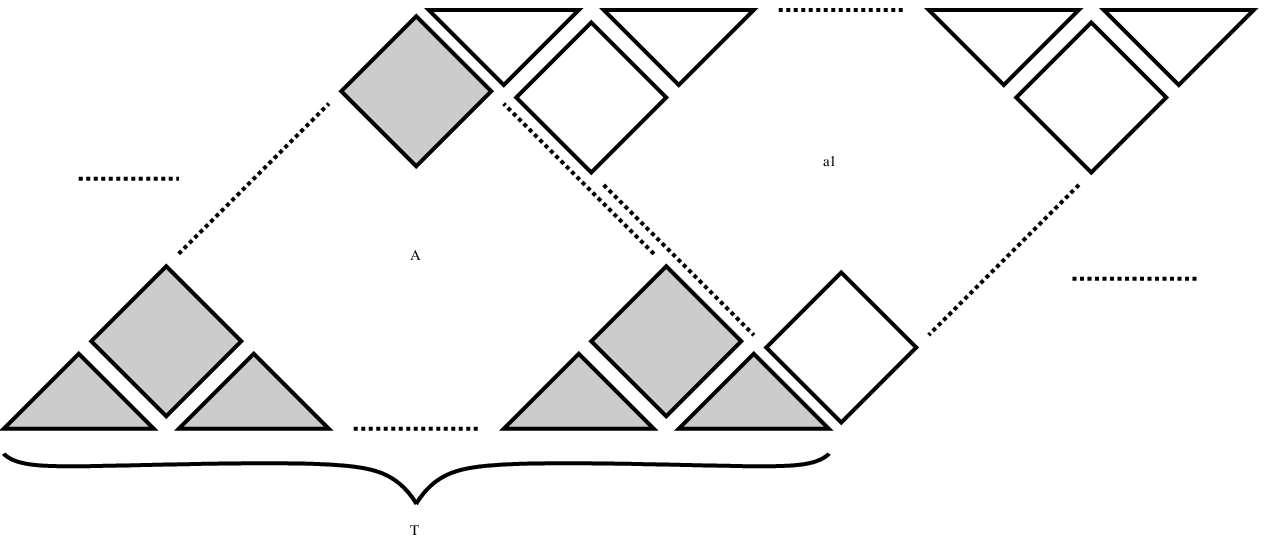}
\end{center}
In this picture the triangles and squares are symbolic representations
for $\bZ A_\infty$ and $\bZ A^\infty_\infty$-components, respectively (see below). 

\medskip

If $\PP=D_\LL$ then the Auslander-Reiten quivers of $\repc(D_\LL)$ and its derived
category have the form:

\vspace*{2mm}

\begin{center}
\psfrag{T}[][]{$\LL$}
\psfrag{A}[][]{$\AA$}
\psfrag{A1}[][]{$\AA[1]$}
\includegraphics[height=3.5cm]{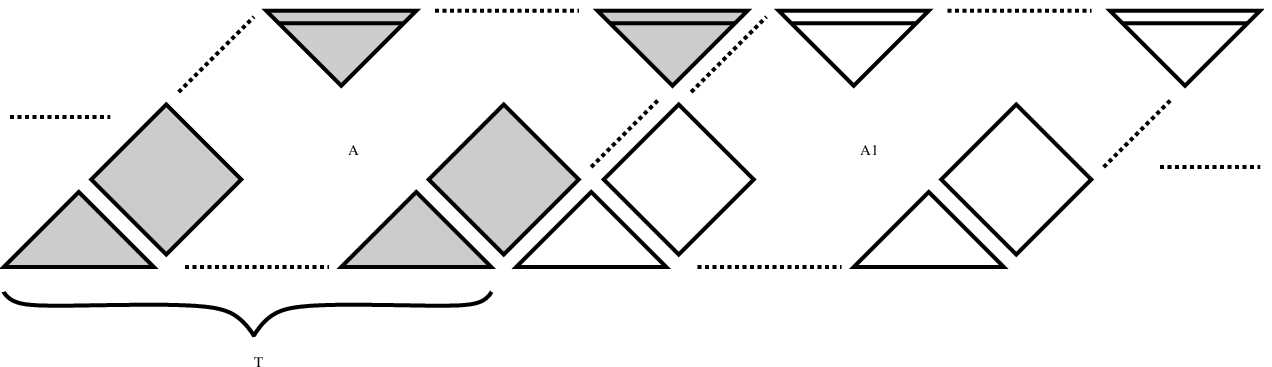}
\end{center}
In this picture the triangles with a double base are symbolic
representation for $\bZ D_\infty$-components (see below).

\medskip

The proof of  Theorem \ref{introtheorem} is quite involved and consists of a number of steps
which we now briefly sketch.

\medskip

\noindent
\textbf{Step 1.}  The following result (see Theorem \ref{theorem:PTS}) is used a various places. Since it does not depend on Serre
duality it  may be of independent interest. 
\begin{theorem} Let $\AA$ be a directed abelian $\Ext$-finite
  hereditary category. Then $\dim \Ext^i_\AA(X,Y)\le 6$ for all
indecomposable  $X,Y\in \AA$ and $i=0,1$. Furthermore if either of the vector spaces $
  \Hom_\AA(X,Y)$ or $\Ext_\AA(Y,X)$ is zero then the other is at most
1-dimensional.
\end{theorem}

\noindent
\textbf{Step 2.}
For the rest of the proof we assume that $\AA$ is a connected directed
hereditary category satisfying Serre duality and we put
$\CC=D^b(\AA)$. Our first aim is to identify the shapes of the
connected components of the Auslander-Reiten quiver of $\CC$. Since such a
component is a stable translation quiver it must be of the form $\bZ
B/G$ where $B$ is an oriented tree and $G$ is acting on $\bZ B$ (\cite{Riedtmann80}). Using an appropriate generalization of the theory of sectional paths (\cite{ARS}) and directedness we deduce:
\begin{enumerate}
\item All components are \emph{standard}, i.e.\ all maps are linear combinations of compositions of
  irreducible ones. In particular all relations can be obtained from the
  mesh relations.
\item $|G|=1$ and furthermore one of the following is true: $B$ is Dynkin,
  $B=A_\infty$, $B=A_{\infty}^{\infty}$ or $B=D_\infty$ (these may be characterized
   as the trees not containing non-Dynkin diagrams).
\end{enumerate}
If $\CC$ has a component $\bZ B$ with $B$ Dynkin then from
connectedness it follows easily that $\CC\cong D^b(\rep(B))$. So below we exclude this case.

\medskip

\noindent
\textbf{Step 3.}  The next step is to understand the maps between
different components. Let $\KK$ be an Auslander-Reiten component of $\CC$.
Since we know all morphisms in $\KK$ we may select a \emph{partial
  tilting set} (\S\ref{section:PTS}) in $\KK$ which generates $\KK$. In this
way we construct a partially ordered set $\PP$ together with an
exact embedding $D^b(\rep(\PP))\r \CC$ whose essential image
contains $\KK$ (and its shifts). The fact that this essential image
usually also contains other components
allows us to obtain information on the interaction between different
components.

\medskip

\noindent
\textbf{Step 4.}
Now we develop the \emph{probing} technique (\S\ref{section:Probing}). Let us say that an indecomposable object $S$ in a $\bZ A_\infty$-component is \emph{quasi-simple} if the middle term of the right Auslander-Reiten triangle built on $S$ is indecomposable.
Using the technique developed in Step 3 we prove that
most indecomposable objects
in $\CC$ have precisely two distinct quasi-simples
mapping to it and these quasi-simples
identify the object uniquely.

\medskip

\noindent
\textbf{Step 5.}
Our next observation is that we if we have a morphism
$X\r Y$ in $\CC$ we can often determine the quasi-simples mapping
to its cone by knowing the quasi-simples mapping to $X$ and $Y$. This gives
us a hold on the triangulated structure of $\CC$.

\medskip

\noindent
\textbf{Step 6. } Finally we use the information gathered in the 
previous steps to construct a tilting set in $\CC$. 
For example if $\CC$ has no $\ZDi$-components then this
tilting set is 
\[
\{X\in \operatorname{ind}(\CC)\mid \Hom_\CC(S,X)\neq 0\}
\]
where $S$ is an arbitrarily chosen quasi-simple in $\CC$. The structure of this tilting set allows us finally to complete the proof of 
Theorem \ref{introtheorem}.

\medskip

\noindent
\textbf{Acknowledgment } I wish to thank Michel Van den Bergh for 
many useful discussions as well as for contributing some ideas.
\section{Notations and conventions}

Throughout this paper, fix an algebraically closed field $k$ of arbitrary characteristic.  All algebras and categories are assumed to be $k$-linear.

If $\mathfrak{a}$ is a small pre-additive category, we denote by $\Mod(\mathfrak{a})$ the abelian category of all left modules, thus all covariant functors from $\mathfrak{a}$ to $\Mod(k)$, the category of vectorspaces over $k$.  The category $\mod(\mathfrak{a})$ is the full subcategory of $\Mod(\mathfrak{a})$ consisting of finitely presented objects.  We will interpret an algebra as a pre-additive category consisting of one object.

For a category $\CC$, we will write $\ind \CC$ for a chosen set of representatives of non-isomorphic indecomposable objects of $\CC$.  If $E$ is an object of $\CC$, then $\ind E$ is a chosen set of representatives of non-isomorphic indecomposable direct summands of $E$.

The ring $\bZ$ is the ring of integers.  We will write $\bZz$ for $\bZ \setminus \{0\}$.

If $\QQ$ is a quiver, a poset, an algebra, or more generally, a pre-additive category, we will write $\QQ^{\circ}$ for the dual quiver, poset, algebra, or pre-additive category.

Given a category $\CC$, the \emph{Auslander-Reiten quiver of $\CC$} is defined as follows: the vertices of $\CC$ are the isomorphism classes  $[X]$ of indecomposable objects $X$.  There is an arrow $[X] \r [Y]$ if there is an irreducible morphism between $X \r Y$.  If $\KK$ is a component of the Auslander-Reiten quiver of $\CC$, also called an \emph{Auslander-Reiten component of $\CC$}, we will often, by abuse of notation, write $X \in \KK$ and $X = Y$, when we mean $[X] \in \KK$ and $[X] = [Y]$, respectively.

For an abelian category or pre-additive category $\AA$, we will write $\Kb \AA$ for the category of bounded complexes modulo homotopy relations.

If $\AA$ is abelian, we write $\Db \AA$ for the bounded derived category of $\AA$.  The category $\Db \AA$ has the structure of a triangulated category.  Whenever we use the word "triangle" we mean "distinguished triangle".

We will say a Krull-Schmidt abelian or triangulated category $\AA$ is \emph{connected} if for all $X, Y \in \ind \AA$ there is a finite, non-oriented path of non-zero morphisms between $X$ and $Y$ or between $X$ and $Y[z]$ for a $z \in \bZ$, respectively.
\section{Preliminaries}

Let $k$ be an algebraically closed field.  A $k$-linear abelian or triangulated category $\AA$ is \emph{Ext-finite} if for all objects $X,Y \in \Ob(\AA)$ one has that $\dim_{k}\Ext^{i}(X,Y)<\infty$ for all $i \in \bN$.  We say that $\AA$ is \emph{hereditary} if $\Ext^{i}(X,Y)=0$ for all $i \geq 2$.

If $\AA$ is an abelian category, we will say that $\AA$ satisfies \emph{Serre duality} if there exists an auto-equivalence $F:\Db\AA \to \Db\AA$, called the \emph{Serre functor}, such that, for all $X,Y\in\Ob(\Db\AA)$, there is an isomorphism
$$\Hom_{\Db\AA}(X,Y) \cong \Hom_{\Db\AA}(Y,FX)^{*}$$
which is natural in $X$ and $Y$ and where $(-)^{*}$ denotes the vector space dual.  Such a functor is necessarily exact (\cite{Bondal89}).

It has been proven in \cite{ReVdB02} that an abelian category $\AA$ has a Serre functor if and only if the category $\Db\AA$ has Auslander-Reiten triangles.  In that case the action of the Serre functor on objects coincides with $\tau[1]$, where $\tau$ is the Auslander-Reiten translation.

We will continue with some general remarks concerning directed categories.  Recall that a category $\CC$ is \emph{directed} if all indecomposable objects are \emph{directing}, thus if for all indecomposable $X\in\Ob(\CC)$, there is no path $X\cong X_{0}\to X_{1} \to X_{2} \to \cdots \to X_{n} \cong X$ of indecomposable objects $X_{i}$, $0\leq i\leq n$, with $\rad(X_{i-1},X_{i})\not=0$ for all $1\leq i\leq n$.  As a result, we may conclude that $\rad(X,X)=0$ and thus $\Hom(X,X)=k$.

Mostly, we will work on the bounded derived category $\Db\AA$ of a directed abelian hereditary category $\AA$.  It is standard that there is a full embedding $\AA\to\Db\AA$ by identifying objects of $\AA$ with complexes of $\Db\AA$ concentrated in degree zero.  We will identify $\AA$ as a subcategory of $\Db\AA$ by this embedding.  In the hereditary case, every element of $\Db\AA$ can be written as a finite direct sum $\oplus_{i} A_{i}[i]$ with $A_{i}\in \Ob\AA$ for all $i$.  Since $\Hom_{\Db\AA}(X,Y[-n])=0$ for $n>0$ and $X,Y\in\Ob(\AA)$, it is easy to check that the category $\Db\AA$ is directed as well.  

We shall formulate certain restrictions on the Hom-sets of directed categories as a lemma.

\begin{lemma}\label{lemma:Directed}
Let $\AA$ be a hereditary directed abelian category.  Consider two indecomposable objects $X,Y\in \Ob\Db\AA$.  If $\Hom(X,Y)\not=0$ then $\Hom(X,Y[z])=0$ and $\Hom(Y,X[z+1])=0$ for all $z\in\bZz$.
\end{lemma}

\begin{proof}
Since $\AA$ is hereditary and $X$ is indecomposable, we know that $X$ is contained in $\AA[z]$ considered as full subcategory of $\Db\AA$, for a certain $z \in \bZ$.  Without loss of generality, assume $X$ to be contained in $\AA[0]$.  Since $\Hom(X,Y) \not= 0$ and $\AA$ is hereditary, we may assume either $Y \in \AA[0]$ or $Y \in \AA[1]$.  If $Y \in \AA[0]$, we have $\Hom(X,Y[z])=0$ for all $z\in\bZz$.  Indeed, the case $z<0$ is clear and the case $z>1$ follows from hereditariness.  Thus assume $\Hom(X,Y[1])\not=0$.  The triangle $\tri{Y}{M}{X}$ built on a nonzero morphism $X \to Y[1]$ yields a path from $Y$ to $X$, contradicting directedness.

We will continue by proving $\Hom(Y,X[z+1])=0$ for all $z\in\bZz$.  This is clear for $z < -1$, follows from directedness for $z=-1$ and from hereditariness when $z > -1$.

Note that in the case $X,Y \in \AA[0]$ we have only used that there exists a path from $X$ to $Y$ to prove that $\Hom(X,Y[z])=0$ for $z\in\bZz$.  The case where $Y\in\AA[1]$ is analogous.  Then $Y[-1]\in \AA$, and the triangle $Y[-1] \to M \to X \to Y$ gives a path from $Y[-1]$ to $X$.  Due to the first part of the proof, this suffices to conclude $\Hom(Y[-1],X[z])=0$ and $\Hom(X,Y[-1][z+1])=0$ or equivalently, $\Hom(Y,X[z+1])=0$ and $\Hom(X,Y[z])=0$.
\end{proof}

Finally, we give a slight reformulation of \cite[Lemma 3]{Ringel05}.

\begin{lemma}\label{lemma:RingelTriangles}
Let $\CC$ be a triangulated category and let $\xymatrix{X\ar[r]^{f}&Y\ar[r]^{g}&Z\ar[r]^{h}&X[1]}$ be a
triangle in $\CC$ with $h\not=0$.  Let $Y=\bigoplus_{i=1}^{n} Y_{i}$ where $Y_{i}$ is not necessarily indecomposable for $i=1,\ldots,n$.  Write $g=(g_{1}, \ldots, g_{n})$ and $f=(f_{1}, \ldots, f_{n})$ with maps $f_{i}:X\to Y_{i}$ and $g_{i}:Y_{i}\to Z$.  Then the following statements are true.
\begin{enumerate}
\item
The morphisms $g_{i}$ are non-invertible for $i=1,\ldots, n$.
\item
If $Z$ is indecomposable, then $f_{i}$ is nonzero for $i=1,\ldots, n$. 
\item
The morphisms $f_{i}$ are non-invertible for $i=1,\ldots, n$.
\item
If $X$ is indecomposable, then $g_{i}$ is nonzero for $i=1,\ldots, n$. 
\end{enumerate}
\end{lemma}

\begin{proof}
\begin{enumerate}
\item
If $g_{i}$ were invertible, then $g$ would be a split epimorphism and $h$ would be zero.
\item
If $f_{i}$ were zero, then consider the following morphisms of triangles
$$
\xymatrix{0\ar[r]\ar[d]&Y_{i}\ar@{=}[r]\ar[d]&Y_{i}\ar[r]\ar[d]^{g_{i}}&0\ar[d] \\
X\ar[r]^{f}\ar[d]&Y\ar[r]^{g}\ar[d]&Z\ar[r]^{h}\ar[d]&X[1]\ar[d] \\
0\ar[r]&Y_{i}\ar@{=}[r]&Y_{i}\ar[r]&0
}$$
Since the two compositions on the left are isomorphisms, so too must the composition $Y_{i}\to Z \to Y_{i}$.  We may conclude that $Y_{i}$ is a  direct summand of $Z$.  But since the split monomorphism $g_{i}:Y_{i} \to Z$ is not an isomorphism (by 1.), this would contradict the indecomposability of $Z$.
\item[(3 \& 4)]
Similar.
\end{enumerate}
\end{proof}

\section{On semi-hereditary pre-additive categories}\label{section:Sh}

For $\frak{a}$  a small pre-additive category we denote by $\Mod(\frak{a})$
the category of left $\frak{a}$-modules. An object $M$ of $\Mod(\frak{a})$
may be represented by a set of objects $M(A)_{A\in \Ob A}$, covariant in $A$.

If $f:\frak{a}\r \frak{b}$ is a functor between small pre-additive
categories then there is an obvious restriction functor
\[
(-)_{\frak{a}}:\Mod(\frak{b})\r \Mod(\frak{a})
\]
which sends $N(B)_{B\in \Ob \frak{b}}$ to $N(f(A))_{A\in \Ob \frak{a}}$. This
restriction functor has a left adjoint 
\[
\frak{b}\otimes_{\frak{a}}-:\Mod(\frak{a})\r \Mod(\frak{b})
\]
which is the right exact functor which sends the projective generators
$\frak{a}(A,-)$ in $\Mod(\frak{a})$ to $\frak{b}(f(A),-)$ in $\Mod(\frak{b})$. As usual if $f$ is fully faithful
we have $(\frak{b}\otimes_{\frak{a}} N)_\frak{a}=N$.

We will also need the following fact.
\begin{lemma} 
\label{thelemma} 
Assume that $f$ is fully faithful and assume that $P$ is a
  summand of an object of the form $\bigoplus_{i=1}^n \frak{b}(B_i,-)$
  with $B_i$ in the essential image of $f$. Then the canonical map
  $\frak{b}\otimes_{\frak{a}} P_{\frak{a}}\r P$ is an isomorphism.
\end{lemma}
\begin{proof} 
$P$ is given by an idempotent $e$ in $\oplus_{i,j} \frak{b}(B_i,B_j)$. Hence
we may write $P$ as the cokernel of 
\[
\bigoplus_{i=1}^n \frak{b}(B_i,-)\xrightarrow{1-e}\bigoplus_{i=1}^n \frak{b}(B_i,-)
\]
The result now follows easily from the right exactness of $\frak{b}\otimes_{\frak{a}}-$.
\end{proof}

Let $M$ be in $\Mod(\frak{a})$. We will say that $M$ is \emph{finitely generated}
if $M$ is a quotient of finitely generated projectives.  Similarly we say that  $M$ is \emph{finitely presented}
if $M$ has a presentation
\[
P\r Q\r M\r 0
\]
where $P,Q$ are finitely generated projectives. It is easy to see that these notions coincides with the ordinary categorical ones.

Dually we will say that $M$ is \emph{cofinitely generated} if it is contained in a cofinitely generated injective. \emph{Cofinitely presented} is defined in a similar way.

The categorical
interpretation of the latter notions is somewhat less clear.  However if $\frak{a}$ is $\Hom$-finite then both finitely and cofinitely presented representations correspond to each other under duality (exchanging $\frak{a}$ and $\frak{a}^\circ$).

We say that a small pre-additive category $\frak{a}$ is \emph{semi-hereditary} if the finitely presented objects $\mod(\frak{a})$ in $\Mod(\frak{a})$ form an abelian and hereditary category.  Following proposition shows that semi-hereditariness is a local property.

\begin{proposition} 
\label{proposition:ShIsLocal}Let $\frak{a}$ be a small pre-additive category such that any
full subcategory of $\frak{a}$ with a finite number of objects is semi-hereditary. Then $\frak{a}$ is itself semi-hereditary.
\end{proposition}
\begin{proof} As usual it is sufficient to prove that the kernel of a
  map between finitely generated projectives $p:P\r Q$ in
  $\Mod(\frak{a})$ is finitely generated projective and splits off.

  Since a finitely generated projective $\frak{a}$-module is a summand of an
  $\frak{a}$-module of the form $\oplus_{i=1}^n \frak{a}(A_i,-)$ we may
 without of loss of generality assume that $p$ is a map of the form
\[
p:\bigoplus_{i=1}^m \frak{a}(A_i,-)\r \bigoplus_{j=1}^n \frak{a}(B_j,-)
\]
Such a map is given by a collection of maps $p_{ji}:B_j\r A_i$.

Let $\frak{b}$ be the full subcategory of $\frak{a}$ containing the
objects $(A_i)_i$, $(B_j)_j$ and let $\FF$ be the filtered
collection of full subcategories of $\frak{a}$ containing $\frak{b}$
and having a finite number of
objects.

 For
$\frak{c}\in \FF$ let $K_{\frak{c}}$ be the
kernel of the map
\[
\bigoplus_{i=1}^m \frak{c}(A_i,-)\r \bigoplus_{j=1}^n \frak{c}(B_j,-)
\]
given by the same $(p_{ji})_{ij}$.
Put $K=\frak{a}\otimes_{\frak{b}} K_{\frak{b}}$. By hypotheses $K_{\frak{b}}$ is finitely
generated and a summand of  $\bigoplus_{i=1}^m\frak{b}(A_i,-)$
and it follows that the analogous facts are true for $K$. So to prove the proposition it is  sufficient to prove that $K$ is the kernel of $p$.

That $K$ is the kernel of $p$ can be checked pointwise.  Hence it is
sufficient to show it for an arbitrary $\frak{c}\in \FF$.  Since
it is easy to see that $(\frak{a}\otimes_{\frak{b}} K_{\frak{b}})_{\frak{c}}=\frak{c}\otimes_{\frak{b}}
K_{\frak{b}}$ we need to show that the canonical map
$\frak{c}\otimes_{\frak{b}} K_{\frak{b}}\r K_{\frak{c}}$ is an isomorphism. Since
$K_{\frak{c}}$ is a summand of $\bigoplus_{i=1}^m \frak{c}(A_i,-)$ and
obviously $K_{\frak{b}}=(K_{\frak{c}})_{\frak{b}}$ this follows from lemma
\ref{thelemma}.
\end{proof}

\begin{remark} Proposition \ref{proposition:ShIsLocal} is false with 
semi-hereditary replaced by hereditary. Furthermore it is not necessarily
true that the functors $\frak{a}\otimes_{\frak{b}}-$ which enter in the proof
are exact. 
\end{remark}
\begin{remark} It is not true that a filtered direct limit of semi-hereditary
rings is semi-hereditary.  A counterexample is given in \cite{Bergman78}.
\end{remark}

Let $\PP$ be a partially ordered set. We associate to $\PP$ 
a pre-additive category $k\PP$ as follows: the object of $k\PP$ are the elements of
$\PP$ and 
\[
(k\PP)(i,j)=
\begin{cases} 
k&\text{if $i\le j$}\\
0&\text{otherwise}
\end{cases}
\]
For $i\le j$ denote the element of $(k\PP)(i,j)$ corresponding to $1\in k$ by $(i,j)$. Composition
of maps in $k\PP$ is defined as  $(j,k)(i,j)=(i,k)$.
We write $\Rep(\PP)$ for $\Mod(k\PP)$ and $\rep(\PP)$ for $\mod(k\PP)$. The objects of $\Rep(\PP)$ are often called $\PP$-representations.

In this article, we will mainly be interested in the category of finitely presented and cofinitely presented representations of a poset $\PP$, which we will denote by $\repc(\PP)$.

Note that if in $\Rep(\PP)$ the finitely generated projectives are cofinitely presented, we have $\repc(\PP) \cong \rep(\PP)$.




We will say that a poset is a \emph{forest} if for all $i,j\in\PP$
such that $i<j$
the \emph{interval}
\[
[i,j]=\{k\in \PP\mid i\le k\le j\}
\]
is totally ordered. It is clear that a subposet of a forest is a forest.

Now, assume that $Q$ is a poset which is a forest. Then any finite subposet $Q_0$ of $Q$ is still a forest and hence $\Rep(Q_0)$ is hereditary. It follows from Proposition \ref{proposition:ShIsLocal} that $\rep(Q)$ and also $\repc(Q)$ are hereditary abelian categories.

We will now proceed to define two posets of special interest.  We will say that a poset is \emph{locally discrete} if no element is an accumulation point.  Thus a linearly ordered poset is locally discrete if and only if for each non-maximal element $i$ there exists an immediate successor $i+1$ and for each non-minimal element $i$ there exists an immediate predecessor $i-1$.  If $\LL$ is a linearly ordered poset, we will denote by $\DL$ the set $\{Q_1 ,Q_2 \} \stackrel{\cdot}{\cup} \LL$ endowed with a poset structure induced by the relation
$$X < Y \asa \left\{\begin{array}{l}
\mbox{$X,Y\in \LL$ and $X <_{\LL} Y$, or}\\
\mbox{$X\in \{Q_1 ,Q_2 \}$ and $Y \in \LL$}
\end{array}\right.
$$

In analogy with the notation used for Dynkin quivers, we will also write $\AL$ for $\LL$.  For the rest of this section, we will assume $\LL$ to be a locally discrete linearly ordered set with no extremal elements, thus not having a maximal nor a minimal element.

\subsection{The category $\repc(\AL)$}\label{section:AL}

These categories have already been considered in \cite{Ringel02}.  In this section we will recall some results.

For all $i,j \in \AL$ with $i \leq j$ we will write 
$$A_{i,j} = \coker((k\AL)(j+1,-) \r (k\AL)(i,-)).$$
It is easily seen that $A_{i,j}$ is cofinitely presented, thus it is an indecomposable object of $\repc(\AL)$.  Following lemma will classify all objects of $\repc(\AL)$.

\begin{lemma}\label{lemma:ModulesAL}
The objects of $\repc(\AL)$ are all isomorphic to finite direct sums of modules of the form $A_{i,j}$.
\end{lemma}

\begin{proof}
We first prove that $A_{i,j}$ is indecomposable.  It is easy to see that the number of indecomposable summands of $A_{i,j}$ is at most $\dim \Hom((k\AL)(i,-),A_{i,j})$.  Applying the functor $(k\AL)(i,-)$ to the exact sequence
$$0 \r (k\AL)(j+1,-) \r (k\AL)(i,-) \r A_{i,j} \r 0$$
yields $\dim \Hom((k\AL)(i,-),A_{i,j}) = 1$.

Conversely, let $X$ be an indecomposable object of $\rep(\AL)$.  Since $X$ is finitely presented in $\Rep(\AL)$, we may choose finitely many projectives generating a full subcategory $\AA$ of $\Rep(\AL)$ containing $X$, such that the embedding $i : \AA \to \Rep(\AL)$ is right exact.  This subcategory $\AA$ is equivalent to $\rep(A_{n})$ for a certain $n \in \bN$, hence
$$X = \coker(f : (k A_{n})(j,-) \r (k A_{n})(i,-))$$
for certain $i,j \in A_{n}$.  We may assume $f\not=0$, since otherwise $X$ would be projective in $\AA$ and in $\repc(\AL)$, and hence will not have a cofinite presentation in $\repc(\AL)$.
\end{proof}

\begin{proposition}\label{proposition:AL}
Let $\LL$ be a locally discrete linearly ordered poset without extremal elements, then the category $\repc(\AL)$ is a connected directed hereditary abelian Ext-finite $k$-linear category satisfying Serre duality.
\end{proposition}

\begin{proof}
It follows from the proof of Proposition \ref{proposition:ShIsLocal} and its dual that $\repc(\AL)$ is abelian and hereditary since it is a full and exact subcategory of $\Rep(\AL)$.

It is easily seen that $\repc(\AL)$ is connected.  We need only check that $\repc(\AL)$ is directed and satisfies Serre duality.

First we show that $\rep(\AL)$ is directed.  Assume there is a cycle $X\cong X_{0}\to X_{1} \to X_{2} \to \cdots \to X_{n} \cong X$ of indecomposable objects $X_{i}$, $0\leq i\leq n$, with $\rad(X_{i-1},X_{i})\not=0$ for all $1\leq i\leq n$.  Since each $X_{i}$ is has a finite presentation in $\Rep(\AL)$, we may choose finitely many projectives generating a full subcategory $\AA$ of $\Rep(\AL)$ containing every $X_{i}$.  Since $\AA$ is equivalent to the directed category $\mod(A_{n})$ for a certain $n\in \bN$, this gives the required contradiction.

Since $\repc(\AL)$ has neither projectives nor injectives (for every indecomposable $A_{i,j}$ there is a non-split epimorphism $A_{i,j+1} \to A_{i,j}$, and a non-slit monomorphism $A_{i,j} \to A_{i-1,j}$), we know by \cite{ReVdB02} that the existence of a Serre functor on $\Db(\repc(\AL))$ is equivalent to the existence of Auslander-Reiten sequences in $\repc(\AL)$.  Let $A_{i,j}$ be an indecomposable object of $\repc(\AL)$.  We claim that the exact sequence
\begin{equation}\label{equation:ARsequence}
0 \r A_{i+1,j+1} \r A_{i+1,j} \oplus A_{i,j+1} \r A_{i,j} \r 0
\end{equation}
is an Auslander-Reiten sequence.  To illustrate this, let $Y$ be an indecomposable object of $\repc(\AL)$ and choose finitely many projectives generating a full subcategory $\AA$ of $\Rep(\AL)$ containing $Y$ and the exact sequence (\ref{equation:ARsequence}).  It is clear that $\AA$ is equivalent to the category $\rep(A_{n})$ for a certain $n\in \bN$ and that the short exact sequence (\ref{equation:ARsequence}) is an almost split exact sequence in $\AA$.  Hence all morphisms $Y \r A_{i,j}$ and $A_{i-1,j-1} \r Y$ factor through the middle term.  This shows that the exact sequence (\ref{equation:ARsequence}) is an Auslander-Reiten sequence in $\repc(\AL)$.
\end{proof}

Finally, we will give the Auslander-Reiten quiver of $\repc(\AL)$.  Therefore, let $\TT$ be any linearly ordered set and consider the poset $\LL = \TXZ$ defined by endowing $\TT \times \bZ$ with the lexicographical ordering.  Thus, for all $t,t' \in \TT$ and $z,z' \in \bZ$, we have
$$(t,z) \leq (t',z') \Leftrightarrow \left\{ \begin{array}{ll}
t < t' \mbox{, or} \\
t=t' \mbox{ and } z \leq z'
\end{array}\right.
$$
It is readily seen that $\LL$ is a locally discrete linearly ordered set with no extremal elements and, conversely, that every such ordered set is constructed in this way.

For every $t \in \TT$ the Auslander-Reiten quiver of $\repc(\AL)$ has a $\ZAi$-component as given in Figure \ref{fig:ZAiOfAL}.  With two distinct elements $t < t'$ correspond a $\ZAii$-component as given in Figure \ref{fig:ZAiiOfAL}.
\begin{figure}
	\centering
$$\xymatrix@C=8pt@R=8pt{\cdots\ar[dr]&&\cdots\ar[dr]&&\cdots\ar[dr]&&\cdots\\
&A_{(t,z),(t,z+2)}\ar[dr]\ar[ur]&&A_{(t,z-1),(t,z+1)}\ar[ur]\ar[dr]&&A_{(t,z-2),(t,z)}\ar[dr]\ar[ur]\\
\cdots\ar[dr]\ar[ur]&&A_{(t,z),(t,z+1)}\ar[dr]\ar[ur]&&A_{(t,z-1),(t,z)}\ar[dr]\ar[ur]&&\cdots \\
&A_{(t,z+1),(t,z+1)}\ar[ur]&&A_{(t,z),(t,z)}\ar[ur]&&A_{(t,z-1),(t,z-1)}\ar[ur]
}$$	
   \caption{A $\ZAi$-component of $\repc(\AL)$}
	\label{fig:ZAiOfAL}
\end{figure}

\begin{figure}
	\centering
$$\xymatrix@C=8pt@R=8pt{\cdots\ar[dr]&&\cdots\ar[dr]&&\cdots\ar[dr]&&\cdots\\
&A_{(t,z),(t',z+2)}\ar[dr]\ar[ur]&&A_{(t,z-1),(t',z+1)}\ar[ur]\ar[dr]&&A_{(t,z-2),(t',z)}\ar[dr]\ar[ur]\\
\cdots\ar[dr]\ar[ur]&&A_{(t,z),(t',z+1)}\ar[dr]\ar[ur]&&A_{(t,z-1),(t',z)}\ar[dr]\ar[ur]&&\cdots \\
&A_{(t,z+1),(t',z+1)}\ar[ur]\ar[dr]&&A_{(t,z),(t',z)}\ar[ur]\ar[dr]&&A_{(t,z-1),(t',z-1)}\ar[ur]\ar[dr]\\
\cdots\ar[ur]&&\cdots\ar[ur]&&\cdots\ar[ur]&&\cdots
}$$	
   \caption{A $\ZAii$-component of $\repc(\AL)$}
	\label{fig:ZAiiOfAL}
\end{figure}
\subsection{The category $\repc(\DL)$}\label{section:DL}

This section closely parallels the previous one, although some arguments are slightly more elaborate.  For all $i,j \in \DL$ with $i \leq j$ we will write
\begin{eqnarray*}
A_{i,j} &=& \coker((k\DL)(j+1,-) \r (k\DL)(i,-)) \\
A^{1}_{j} &=& \coker((k\DL)(j+1,-) \r (k\DL)(Q_{1},-)) \\
A^{2}_{j} &=& \coker((k\DL)(j+1,-) \r (k\DL)(Q_{2},-)) \\
B_{i,j} &=& \coker((k\DL)(j+1,-) \oplus (k\DL)(i+1,-) \r (k\DL)(Q_{1},-) \oplus (k\DL)(Q_{2},-))
\end{eqnarray*}
where in the definition of $B_{i,j}$ we assume $i \not= j$.

It is easy to see that $A_{i,j}, A^{1}_{j}, A^{2}_{j}$ and $B_{i,j}$ are also cofinitely presented, hence they are objects of $\repc(\DL)$.  In following lemma we prove that these are all indecomposable objects.

\begin{lemma}
The objects of $\repc(\DL)$ are all isomorphic to finite direct sums of modules of the form $A_{i,j}, A^{1}_{j}, A^{2}_{j}$ or $B_{i,j}$.
\end{lemma}

\begin{proof}
Analogue to the proof of Lemma \ref{lemma:ModulesAL}
\end{proof}

\begin{proposition}
Let $\LL$ be a locally discrete linearly ordered poset without extremal elements, then the category $\repc(\DL)$ is a connected directed hereditary abelian Ext-finite category satisfying Serre duality.
\end{proposition}

\begin{proof}
Analogue to the proof of Proposition \ref{proposition:AL}.  
\end{proof}

As in the case of $\repc(\AL)$, we will give a description of the Auslander-Reiten components of $\repc(\DL)$.  Let $\LL = \TXZ$, for a certain linearly ordered poset $\TT$.  For every $t \in \TT$ the Auslander-Reiten quiver of $\repc(\DL)$ has a $\ZAi$-component and $\ZDi$-component as given in Figures \ref{fig:ZAiOfDL} and \ref{fig:ZDiOfDL}.  With two distinct elements $t < t'$ correspond two $\ZAii$-component as given in Figures \ref{fig:ZAii1OfDL} and \ref{fig:ZAii2OfDL}.

\begin{figure}
	\centering
$$\xymatrix@C=8pt@R=8pt{\cdots\ar[dr]&&\cdots\ar[dr]&&\cdots\ar[dr]&&\cdots\\
&A_{(t,z),(t,z+2)}\ar[dr]\ar[ur]&&A_{(t,z-1),(t,z+1)}\ar[ur]\ar[dr]&&A_{(t,z-2),(t,z)}\ar[dr]\ar[ur]\\
\cdots\ar[dr]\ar[ur]&&A_{(t,z),(t,z+1)}\ar[dr]\ar[ur]&&A_{(t,z-1),(t,z)}\ar[dr]\ar[ur]&&\cdots \\
&A_{(t,z+1),(t,z+1)}\ar[ur]&&A_{(t,z),(t,z)}\ar[ur]&&A_{(t,z-1),(t,z-1)}\ar[ur]
}$$
   \caption{A $\ZAi$-component of $\repc(\DL)$}
	\label{fig:ZAiOfDL}
\end{figure}

\begin{figure}
	\centering
$$\xymatrix@C=8pt@R=8pt{
\cdots\ar[rdd] && A^{1}_{(t,i+1)}\ar[rdd] && A^{1}_{(t,i)}\ar[rdd] && \cdots \\
\cdots\ar[rd]  && A^{2}_{(t,i+1)}\ar[rd] && A^{2}_{(t,i)}\ar[rd] && \cdots \\
&B_{(t,i+1),(t,i+2)}\ar[ruu]\ar[ru]\ar[rd]&&B_{(t,i),(t,i+1)}\ar[ruu]\ar[ru]\ar[rd]&&B_{(t,i-1),(t,i)}\ar[ruu]\ar[ru]\ar[rd]\\
\cdots\ar[ru]\ar[rd]&&B_{(t,i),(t,i+2)}\ar[ru]\ar[rd]&&B_{(t,i-1),(t,i+1)}\ar[ru]\ar[rd]&&\cdots\\
&B_{(t,i),(t,i+3)}\ar[ru]\ar[rd]&&B_{(t,i-1),(t,i+2)}\ar[ru]\ar[rd]&&B_{(t,i-2),(t,i+1)}\ar[ru]\ar[rd]\\
\cdots\ar[ru]&&\cdots\ar[ru]&&\cdots\ar[ru]&&\cdots
}$$
   \caption{A $\ZDi$-component of $\repc(\DL)$}
	\label{fig:ZDiOfDL}
\end{figure}

\begin{figure}
	\centering
$$\xymatrix@C=8pt@R=8pt{\cdots\ar[dr]&&\cdots\ar[dr]&&\cdots\ar[dr]&&\cdots\\
&A_{(t,z),(t',z+2)}\ar[dr]\ar[ur]&&A_{(t,z-1),(t',z+1)}\ar[ur]\ar[dr]&&A_{(t,z-2),(t',z)}\ar[dr]\ar[ur]\\
\cdots\ar[dr]\ar[ur]&&A_{(t,z),(t',z+1)}\ar[dr]\ar[ur]&&A_{(t,z-1),(t',z)}\ar[dr]\ar[ur]&&\cdots \\
&A_{(t,z+1),(t',z+1)}\ar[ur]\ar[dr]&&A_{(t,z),(t',z)}\ar[ur]\ar[dr]&&A_{(t,z-1),(t',z-1)}\ar[ur]\ar[dr]\\
\cdots\ar[ur]&&\cdots\ar[ur]&&\cdots\ar[ur]&&\cdots
}$$	
   \caption{The first $\ZAii$-component of $\repc(\DL)$}
	\label{fig:ZAii1OfDL}
\end{figure}

\begin{figure}
	\centering
$$\xymatrix@C=8pt@R=8pt{\cdots\ar[dr]&&\cdots\ar[dr]&&\cdots\ar[dr]&&\cdots\\
&B_{(t,z),(t',z+2)}\ar[dr]\ar[ur]&&B_{(t,z-1),(t',z+1)}\ar[ur]\ar[dr]&&B_{(t,z-2),(t',z)}\ar[dr]\ar[ur]\\
\cdots\ar[dr]\ar[ur]&&B_{(t,z),(t',z+1)}\ar[dr]\ar[ur]&&B_{(t,z-1),(t',z)}\ar[dr]\ar[ur]&&\cdots \\
&B_{(t,z+1),(t',z+1)}\ar[ur]\ar[dr]&&B_{(t,z),(t',z)}\ar[ur]\ar[dr]&&B_{(t,z-1),(t',z-1)}\ar[ur]\ar[dr]\\
\cdots\ar[ur]&&\cdots\ar[ur]&&\cdots\ar[ur]&&\cdots
}$$	
   \caption{The second $\ZAii$-component of $\repc(\DL)$}
	\label{fig:ZAii2OfDL}
\end{figure}
\section{Partial Tilting Sets}\label{section:PTS}

In this section we shall assume $\AA$ is a $k$-linear abelian Ext-finite category, not necessarily satisfying Serre duality.  We will say that the set $\{P_{i}\}_{i\in I} \subseteq \ind\Db\AA$ is a \emph{partial tilting set} if $\Hom(P_{i},P_{j}[z])=0$, for all $z\in\bZz$ and all $i,j\in I$. 

In the rest of this article, we will often use next theorem.  Recall that a category is called \emph{Karoubian} if the category has finite direct sums and idempotents split.  A small pre-additive category $\frak{a}$ is \emph{coherent} if the finitely presented objects $\mod(\frak{a})$ in $\Mod(\frak{a})$ form an abelian category.

\begin{theorem}\label{theorem:PTS}
Let $\AA$ be a $k$-linear abelian category, $\{P_{i}\}_{i\in I}$ a partial tilting set of $\Db\AA$ and $\frak{a}$ the pre-additive category given by $\{P_{i}\}_{i\in I}$ as a full subcategory of $\Db\AA$.  Assume that $\frak{a}^\circ$ is Karoubian and coherent, and that $\mod(\frak{a}^\circ)$ has finite global dimension, then there is a full exact embedding $\Db(\mod{\frak{a}^{\circ}})\to\Db\AA$ sending $\Hom(P_{i},-)$ to $P_{i}$.
\end{theorem}

\begin{proof}
Due to the conditions on the pre-additive category $\frak{a}$, we know that $\Db(\mod{\frak{a}^\circ})$ is equivalent to $\Kb\frak{a}$.  It thus suffices to construct a full and exact embedding $\Kb \frak{a}\to\Db\AA$.  It is well known that the category $\Ind\AA$ of left exact contravariant functors from $\AA$ to $\Mod{k}$ is a $k$-linear Grothendieck category and that the Yoneda embedding of $\AA$ into $\Ind\AA$ is a full and exact embedding.  By \cite[Prop. 2.14]{LoVdB04}, this embedding extends to a full and exact embedding $\Db\AA\to\Db\Ind\AA$. As a Grothendieck category, $\Ind\AA$ has enough injectives and we may, by \cite[Prop. 10.1]{Rickard89}, consider the full and exact embedding $\Kb \frak{a} \to \Db\Ind\AA$ that extends the embedding of $\frak{a}$ in $\Db\Ind\AA$.  Induction over triangles shows that the essential image of $\Kb \frak{a}$ lies in $\Db\AA$.
\end{proof}

Next, let $\AA$ be a directed $k$-linear abelian hereditary Ext-finite category, not necessarily satisfying Serre duality.  Given a pair of indecomposable objects, $X$ and $Y$, in $\Db\AA$ we will use Theorem \ref{theorem:PTS} to find a full and exact subcategory $\BB$ of $\Db\AA$, such that $X,Y \in \BB$ and $\BB \cong \Db(\mod A)$ for a certain finite dimensional $k$-algebra $A$.

\begin{lemma}\label{lemma:CanonicalMorphism}
Let $\AA$ be a directed Ext-finite hereditary $k$-linear category and let $X,Y\in\ind\Db\AA$ such that $\Hom(X,Y)\not=0$.  Consider the triangle 
\begin{equation}\label{equation:Canonical}
\tri{E}{X\otimes Hom(X,Y)}{Y}
\end{equation}
built on the canonical map $X \otimes \Hom(X,Y) \to Y$, then $(\ind E) \cup \{X\}$ is a partial tilting set.
\end{lemma}

\begin{proof}
In order to ease notation, write $\Hom(X,Y)=V$.  We will prove that $(\ind E) \cup \{X\}$ is a partial tilting set by applying Hom-functors to triangle (\ref{equation:Canonical}).  Out of the long exact sequence given by $\Hom(X,-)$ and directedness we deduce that $\Hom(X,E[z])=0$ for all $z\not=0,1$.  For $z=1$, consider the following exact sequence
$$\Hom(X,X \otimes V) \to \Hom(X,Y) \to \Hom(X,E[1]) \to \Hom(X,(X \otimes V)[1]).$$
Since the map $\Hom(X,X \otimes V) \to \Hom(X,Y)$ is an isomorphism and since due to directedness $\Hom(X,(X \otimes V)[1])=0$, we have $\Hom(X,E[1])=0$.  Thus, $\Hom(X,E[z])=0$, for all $z \in \bZz$.

Next, the functor $\Hom(-,X)$ yields a long exact sequence from which we easily deduce that $\Hom(E,X[z])=0$, for all $z \in \bZz$.

Finally, we are left to prove that $\Hom(E,E[z])=0$, for all $z \in \bZz$.  To this end, we first apply the functor $\Hom(Y,-)$ to triangle (\ref{equation:Canonical}).  Using Lemma \ref{lemma:Directed} to see that $\Hom(Y,Y[z])=0$ and $\Hom(Y,X[z+1])=0$ for all $z \in \bZz$, we may deduce that $\Hom(Y,E[z+1])=0$ for all $z \in \bZz$.  One may now readily see that the long exact sequence one acquires by applying $\Hom(-,E)$ to triangle (\ref{equation:Canonical}) yields $\Hom(E,E[z])=0$ for all $z \in \bZz$.  This proves the assertion.
\end{proof}

\begin{theorem}\label{theorem:BoundedHomDim}
Let $\AA$ be a directed abelian Ext-finite $k$-linear hereditary category and let $X$ and $Y$ be two indecomposable objects of $\AA$, then $\dim\Hom(X,Y)\leq 6$ and $\dim\Ext(X,Y)\leq 6$.  If $\Ext(Y,X)=0$, respectively $\Hom(Y,X)=0$, then $\dim\Hom(X,Y)\leq 1$, respectively $\dim\Ext(X,Y)\leq 1$.
\end{theorem}

\begin{proof}
We will work on the derived category $\Db\AA$.  Possibly by renaming $Y[1]$ to $Y$, it suffices to prove that $\dim\Hom_{\Db\AA}(X,Y)\leq 6$ and $\dim\Hom_{\Db\AA}(X,Y)\leq 1$ if $\Hom_{\Db\AA}(Y,X[1])=0$.

We may assume $\Hom(X,Y)\not=0$.  Lemma \ref{lemma:Directed} then yields that $\Hom(X,Y[z])=0$ and $\Hom(Y,X[z+1])=0$ for all $z \in \bZz$.  If furthermore $\Hom(Y,X[1])=0$, then $\{X,Y\}$ is a partial tilting set and, due to Theorem \ref{theorem:PTS}, we know $A = [\End(X\oplus Y)]^{\circ}$ is a representation-directed algebra, i.e.\ $\mod A$ is a directed category.  From this we may deduce that $\dim\Hom(X,Y)=1$.

If $\Hom(Y,X[1])\not=0$, then we turn our attention to the triangle $\tri{E}{X\otimes \Hom(X,Y)}{Y}$.  Lemma \ref{lemma:CanonicalMorphism} yields that $(\ind E) \cup \{X\}$ is a partial tilting set.  Denote the algebra $[\End (E \oplus X)]^{\circ}$ by $A$.  Theorem \ref{theorem:PTS} then gives a full and exact embedding $i:\Db \mod A \to \Db \AA$.  This shows that $A$ is a representation-directed algebra.  Let $P$ and $Q$ be the projective objects of $\mod A$ corresponding to $E$ and $X$, respectively, under $i$.  Since $i$ is exact, we know $R = \cone(P \to Q \otimes \Hom(X,Y))$ corresponds to $Y$ under $i$.

By applying the functor $\Hom(-,Y)$ to the triangle $\tri{E}{X\otimes \Hom(X,Y)}{Y}$ one sees that $\Hom(E\oplus X,Y[z])=0$ for all $z \in \bZz$.  Indeed, by Lemma \ref{lemma:Directed} we have $\Hom(X,Y[z])=0$ and then the long exact sequence yields $\Hom(E,Y[z])=0$.

Since $\Hom(E\oplus X,Y[z])=0$ for all $z \in \bZz$, we deduce that $\Hom(P\oplus Q,R[z])=0$ for all $z \in \bZz$ and hence we may interpret $R$ as an $A$-module.  Thus $\dim\Hom(Q,R)$ is the number of times the top of $Q$ occurs in the Jordan-H\"{o}lder decomposition of $R$.  This, and thus also $\dim\Hom(X,Y)$ is bounded by 6, by \cite[2.4 (9'')]{Ringel84}
\end{proof}
\section{About the Auslander-Reiten components}

In this section, let $\AA$ be a connected directed hereditary abelian $k$-linear Ext-finite category satisfying Serre duality, and write $\CC=\Db\AA$.  Since $\AA$ satisfies Serre duality, the bounded derived category $\CC=\Db\AA$ is stable.  Due to \cite{Riedtmann80}, we may state that the only possible components of the Auslander-Reiten quivers are of the form $\bZ B/G$, where $G$ is an admissible subgroup of $\Aut(\bZ B)$ and $B$ is an oriented tree.  We wish to show that, in our context, the only possible components are $\ZAi$, $\ZAii$, $\ZDi$, or $\bZ Q$ where $Q$ is Dynkin, and that all such components are standard, i.e.\ $\infrad(X,Y)=0$ when $X$ and $Y$ are contained in the same Auslander-Reiten component.  We will need some preliminary results, starting with a definition from \cite{ARS}.

\begin{definition}
Let $I$ be the integers in one of the intervals $]-\infty,n]$, $[n,+\infty[$, $[m,n]$ for $m<n$ or $\{1,\ldots,n\}$ modulo $n$.  Let $\cdots\to A_{i}\to A_{i+1}\to\cdots$ be a sequence of irreducible morphisms between indecomposables with each index in $I$.  The sequence is said to be sectional if $\tau A_{i+2}\not\cong A_{i}$ whenever both $i$ and $i+2$ are in $I$.  The corresponding path in the Auslander-Reiten quiver is said to be a sectional path.
\end{definition}

\begin{lemma}\label{lemma:noExtsOnSectionalPaths}
Let $A$ and $B$ be indecomposables of $\CC$.  If there exists a sectional path $A = A_{0}\to\cdots\to A_{n} = B$, then $\dim\Hom(A,B)=1$.
\end{lemma}

\begin{proof}
We start by noting that $\Hom(A,B)\not=0$ since the composition of irreducible maps from the same sectional sequence is non-zero (\cite[Proposition 3.2]{XiZh}).  Next, we will prove by induction on $n$ that $\Hom(A\oplus A_{n},(A \oplus A_{n})[z])=0$ for all $z\in\bZz$ and thus, by Lemma \ref{theorem:BoundedHomDim}, that $\dim\Hom(A,B)=1$.  

First, assume $n = 1$, or equivalently, $A_{1} = B$.  There is an irreducible morphism $A\to B$ and hence also a morphism $\tau B\to A$.  Therefore, by directedness, $\Hom(A,\tau B)=0$ and by Serre duality, $\Hom(B,A[1])=0$.  Also, since $\Hom(A,B)\not=0$, Lemma \ref{lemma:Directed} yields $\Hom(A,B[z])=0$ and $\Hom(B,A[z+1])=0$ for all $z\in\bZz$.  Combining those two facts gives $\Hom(A\oplus B, (A\oplus B)[z])=0$ for all $z\in\bZz$.  Applying Theorem \ref{theorem:BoundedHomDim} then yields the assertion in case $n=1$.  

Next, assume the assertion has been proven for $n\in\{1,2,\ldots,k-1\}$, we wish to prove the case $n=k$.  First of all note that since $\Hom(A\oplus A_{k-1},(A \oplus A_{k-1})[z])=0$ for all $z \in \bZz$ we have, using Serre duality, $\dim\Hom(A_{k-1},A[1]) = \dim\Hom(A, \tau A_{k-1})= 0$.  Considering the almost split triangle  $\tau A_{k-1} \to M \to A_{k-1} \to \tau A_{k-1}[1]$, and using the fact that $\dim\Hom(A,A_{k-1}) = 1$, one has $\dim\Hom(A,M)=1$.  Since $M$ has both $A_{k-2}$ and $\tau A_{k}$ as directs summands and it has already been proven that $\dim\Hom(A,A_{k-2})=1$, it follows easily that $\Hom(A,\tau A_{k})=0$, and thus $\Hom(A_{k},A[1])=0$.  It has already been noted that $\Hom(A,A_{k})\not=0$ since they lie in the same sectional sequence and thus, by Lemma \ref{lemma:Directed}, that $\Hom(A_{k},A[z])=0$ and $\Hom(A,A_{k}[z+1])=0$ for all $z\in\bZz$.  Combining this with the earlier proven $\Hom(A_{k},A[1])=0$, we arrive at $\Hom(A\oplus A_{k},(A \oplus A_{k})[z])=0$ for all $z\in\bZz$.  Finally, we may invoke Theorem \ref{theorem:BoundedHomDim} to conclude that $\dim\Hom(A,A_{k})=1$.
\end{proof}

\begin{proposition}\label{proposition:Sectional}
Let $A$ and $B$ be indecomposable objects of $\CC$.  If there exists a sectional path $\xymatrix@1{A = A_{0}\ar[r]^-{\alpha_{0}} & A_{1}\ar[r]^-{\alpha_{1}} & \cdots\ar[r]^-{\alpha_{n-1}} & A_{n} = B}$, then $\dim \Hom(A,B)=1$ and $\infrad(A,B)=0$.  Also, there is only one sectional path from $A$ to $B$ in the Auslander-Reiten quiver.
\end{proposition}

\begin{proof}
It has already been proven in Lemma \ref{lemma:noExtsOnSectionalPaths} that $\dim \Hom(A,B)=1$.
Now, if $\infrad(A,B)\not=0$, then $\dim\infrad(A,B)=1$.  Let $\tau B\to M \to B \to \tau B[1]$ be the Auslander-Reiten triangle built on $B$.  A non-zero morphism $f\in\infrad(A,B)$ should factor through $M$ or, more precisely, through the direct summand $A_{n-1}$ of $M$.  Indeed, since Lemma \ref{lemma:noExtsOnSectionalPaths} implies $\Hom(A,\tau B)=0$ we easily obtain $\dim\Hom(A,M)=1$ and since $\dim\Hom(A,A_{n-1})=1$ we may write $f=\alpha_{n-1} \circ f_{n-1}$, with $f_{n-1}\in\infrad(A,A_{n-1})$.  By iteration, one has that $f=\alpha_{n-1}\circ\ldots\circ\alpha_{1}\circ f_{1}$, with $f_{1}\in\infrad(A,A_{1})$, clearly a contradiction.  Thus $\infrad(A,B)=0$.

Finally, we will show that there can be at most one sectional path from $A$ to $B$ in the Auslander-Reiten quiver.  Consider there to be two sectional paths, $A\to A_{1}\to\ldots\to A_{n-1}\to B$ and $A\to A'_{1}\to\ldots\to A'_{n-1}\to B$.  We may assume that $A_{n-1}$ and $A'_{n-1}$ are not isomorphic.  Consider the Auslander-Reiten triangle $\tri{\tau B}{M}{B}$.  Applying the functor $\Hom(A,-)$ we get the exact sequence 
$$\Hom(A,\tau B)\to\Hom(A,M)\to\Hom(A,B)\to 0$$
Now, both $A_{n-1}$ and $A'_{n-1}$ are direct summands of $M$, thus $\dim\Hom(A,M)\geq 2$.  But Lemma \ref{lemma:noExtsOnSectionalPaths} yields $\dim\Hom(A,B)=1$ and $\dim\Hom(A,\tau B)=\dim\Hom(B,A[1])=0$, clearly a contradiction.
\end{proof}

We will now discuss the form of the components that can occur in the Auslander-Reiten quiver of the category $\CC$.  Recall from \cite{Riedtmann80} that a stable component $\KK$ from the Auslander-Reiten quiver of $\CC$ is covered by $\pi:\bZ B \to \KK$, where $B$ is defined as follows : fix a vertex $X$ from $\KK$, then the vertices of $B$ are defined to be all (finite, non-trivial) sectional paths of $\KK$ starting at $X$, and there is an arrow in $B$ from the sectional path $X\to\cdots\to Y$ to the sectional path $X\to\cdots Y\to Z$.  With these definitions, it is clear that $B$ is a tree with a unique source.  There also is a morphism $f: B\to\KK$ by mapping a sectional path $X\to\cdots\to Y$ to $Y$.  This morphism $f$ extends to the covering $\pi:\bZ B \to \KK$ of translation quivers given by
$$(z, X\to\cdots\to Y) \mapsto \tau^{-z} Y.$$
In the following lemma, we will prove that the map $\pi$ is injective, such that $\bZ B \cong \KK$.

\begin{lemma}\label{lemma:ZB}
Every component of the Auslander-Reiten quiver of $\CC$ is isomorphic to $\bZ B$, as stable translation quivers, for a certain oriented tree $B$ with a unique source.
\end{lemma}

\begin{proof}
As stated before, we need only to prove that the map 
$$\pi:\bZ B \to \KK : (z, X\to\cdots\to Y) \mapsto \tau^{-z} Y$$
is injective.  Consider $(z, X\to\cdots\to Y), (z', X\to\cdots\to Y') \in \bZ B$.  Seeking a contradiction, assume that $(z, X\to\cdots\to Y) \not= (z', X\to\cdots\to Y')$ and $\pi (z, X\to\cdots\to Y) = \pi (z', X\to\cdots\to Y')$, thus $\tau^{-z} Y = \tau^{-z'} Y'$.  Thus we assume there to be in $\KK$ two sectional paths starting in the same vertex, and ending in the same $\tau$-orbit.

We will consider the sectional paths
$$X=A_{0} \to A_{1} \to \cdots \to A_{n-1} \to A_{n} = Y$$
and
$$X=B_{0} \to B_{1} \to \cdots \to B_{m-1} \to B_{m} = Y'=\tau^{z'-z}Y$$
where we may without loss of generality assume that $A_{i} \not= \tau^{k} B_{j}$ for $1 \leq i \leq n-1$, for $1 \leq j \leq m-1$ and for all $k \in \bZ$.

We will consider two separate cases.  First, assume that $z'-z \geq n$.  In that case, we have a path from $\tau^{z'-z}Y$ to $\tau^{n}Y$ and a path
$$\tau^{n}Y \to \tau^{n-1} A_{n-1} \to \tau^{n-2} A_{n-2}\to \cdots \to A_{0}= B_{0} \to B_{m} = Y' = \tau^{z'-z} Y$$
contradicting directedness.

If $z'-z < n$, then we find two different sectional paths
$$X = B_{0} \to B_{1} \to \cdots \to B_{m} = \tau^{z'-z} Y \to \tau^{z'-z-1} A_{n-1} \to \cdots \to \tau A_{n-(z'-z)+1} \to A_{n-(z'-z)}$$
and
$$X = A_{0} \to A_{1} \to \cdots \to A_{n-(z'-z)}$$
from $X$ to $A_{n-(z-z')}$ contradicting Proposition \ref{proposition:Sectional}.
\end{proof}

We are now ready to prove the main theorem of this section.

\begin{theorem}\label{theorem:Components}
Let $\AA$ be a directed hereditary abelian $k$-linear Ext-finite category satisfying Serre duality.  Each component of the Auslander-Reiten quiver of $\CC=\Db\AA$ is both 
\begin{enumerate}
\item standard, and
\item of the form $\ZAi$, $\ZAii$, $\ZDi$ or $\bZ Q$, where $Q$ is a quiver of Dynkin type.
\end{enumerate}
\end{theorem}

We will split the proof of this theorem over the next two lemmas.  In the next lemma, we will denote by $d(a,b)$ the usual graph-theoretical distance between vertices $a$ and $b$.

\begin{lemma}\label{lemma:ComponentStandard}
Each component of the Auslander-Reiten quiver of $\CC$ is standard.
\end{lemma}

\begin{proof}
Let $X$ and $Y$ be two indecomposable objects in an Auslander-Reiten component $\KK$ of $\CC$.  In order to prove $\infrad(X,Y)=0$ we will write $\KK$ as $\bZ B$ such that $X$ corresponds to $(0,b)$ where $b$ is the source of the tree $B$.  We will consider two cases.

The first case is where $Y$ has coordinates $(n,v_{Y})$ with $n \geq 0$ and $v_{Y} \in B$.  If $n=0$ then $\infrad(X,Y)=0$ as a consequence of Proposition \ref{proposition:Sectional}.  Now assume $n > 0$.  If $\infrad(X,Y) \not= 0$, then consider the Auslander-Reiten triangle $\tri{\tau Y}{M}{Y}$.  There is at least one indecomposable summand $Y'$ of $M$ such that $\infrad(X,Y') \not= 0$.  The coordinates of $Y'$ are either $(n,v_{Y'})$ where $d(b,v_{Y'})=d(b,v_{Y})-1$ or $(n-1,v_{Y'})$ where $d(b,v_{Y'})=d(b,v_{Y})+1$.  Since $d(b,v_{Y})$ is finite, iteration gives $\infrad(X,Z) \not= 0$ for a certain $Z$ with coordinates $(0,v_{Z})$. This contradicts Proposition \ref{proposition:Sectional}.

The last case is where $Y$ has coordinates $(-n,v_{Y})$ with $n > 0$.  We will proceed by induction on $n$ to prove that $\Hom(X,Y)=0$.  First, we consider $n=1$, thus let $Y$ have coordinates $(-1,v_{Y})$.  We will use a second induction argument on $d=d(b,v_{Y})$.  If $d=0$, then $b=v_{Y}$ and there is a path from $Y$ to $X$.  Directedness then implies $\Hom(X,Y)=0$.  Now, assume $d \geq 1$.  Choose a direct predecessor $v_{Z}$ of $v_{Y}$.  Let $Z \in \KK$ be the indecomposable object corresponding to the coordinates $(-1,v_{Z})$, then $Y$ is a direct summand of $M$ where $M$ is defined by the Auslander-Reiten triangle $\tri{Z}{M}{\tau^{-1}Z}$.

By the induction on $d$ above we know that $\Hom(X,Z)=0$ and by Proposition \ref{proposition:Sectional} that $\dim \Hom(X,\tau^{-1}Z) = 1$ so we infer that $\dim \Hom(X,M) = 1$ if $X \not\cong \tau^{-1} Z$, and $\dim \Hom(X,M) = 0$ if $X \cong \tau^{-1} Z$.  Thus there is at most one indecomposable direct summand $X'$ of $M$ such that $\Hom(X,X')\not= 0$.  This needs to be an indecomposable object lying on a sectional path from $X$ to $\tau^{-1}Z$, hence $X' \not\cong Y$ since $Y$ has coordinates $(-1,v_{Y})$.  We conclude $\Hom(X,Y)=0$.

We will now proceed with the induction on $n$, thus let $Y$ have coordinates $(-n,v_{Y})$ with $n > 1$.  Also in this case, we will use a second induction argument on $d=d(b,v_{Y})$.  As above, if $d=0$, then $b=v_{Y}$ and there is a path from $Y$ to $X$.  Directedness then implies $\Hom(X,Y)=0$.  Thus we may assume $d \geq 1$.  Let $v_Z$ be a direct predecessor of $v_Y$ in the tree $B$ and let $Z \in \KK$ be the indecomposable object corresponding to the coordinates $(-n,v_{Z})$.  Again, $Y$ is a direct summand of $M$ where $M$ is defined by the Auslander-Reiten triangle $\tri{Z}{M}{\tau^{-1}Z}$.

By induction on $d$ and $n$ we know $\Hom(X,Z)=0$ and $\Hom(X,\tau^{-1}Z)=0$, respectively, and hence also $\Hom(X,M)=0$.  We conclude $\Hom(X,Y)=0$.
\end{proof}

\begin{lemma}
Each component of the Auslander-Reiten quiver of $\CC$ is of the form $\ZAi$, $\ZAii$, $\ZDi$ or $\bZ Q$, where $Q$ is a quiver of Dynkin type.
\end{lemma}

\begin{proof}
This is an immediate consequence of Proposition \ref{proposition:SectionPTS}.
\end{proof}

For the statement of the next proposition, we will say that a subquiver $Q$ of a stable translation quiver $T$ is a {\em section} if the embedding of quivers $Q \to \TT$ lifts to an isomorphism $\bZ Q \cong T$ as stable translation quivers, hence $Q$ contains exactly one object from every $\tau$-orbit of $\TT$ and if $X \in Q$ and there is an arrow $X \to Y$ in $\TT$, then either $Y \in Q$ or $\tau Y \in Q$.

\begin{proposition}\label{proposition:SectionPTS}
Let $\KK$ be a component of the Auslander-Reiten quiver of $\CC$, and let $Q$ be a section of $\KK$.  Then the vertices of $Q$ form a partial tilting set.
\end{proposition}

\begin{proof}
Let $X$ and $Y$ be two indecomposable objects from the section $Q$.  We must show that $\Hom(X,Y[z])=0$ for $z \in \bZz$.  Without loss of generality, assume $X \in \Ob\AA[0]$.  It is clear that there are $i,j \in \bN$ such that there is a sectional path from $\tau^{i} Y$ to $X$ and from $X$ to $\tau^{-j} Y$ and thus, by Proposition \ref{proposition:Sectional} and the fact that $\AA$ is hereditary, we may conclude $\tau^{i} Y \in \Ob\AA[-1]$ or $\tau^{i} Y \in \Ob\AA[0]$, and $\tau^{-j} Y \in \Ob\AA[0]$ or $\tau^{-j} Y \in \Ob\AA[1]$.  Since there are paths from $\tau^{-j} Y$ to $Y$ and from $Y$ to $\tau^{i} Y$, we may infer that $Y \in \AA[-1]$, $Y \in \AA[0]$, or $Y \in \AA[1]$.  Hence $\Hom(X,Y[z])=0$ for $z < -1$ and $z > 2$.  We will proceed to show that $\Hom(X,Y[z])=0$ for $z \in \{-1,1,2\}$.

Therefore, we will show there exists an indecomposable $Z \in \KK$ such that there are sectional paths from $Z$ to both $X$ and $Y$.  Indeed, let $n \in \bN$ be the smallest natural number such that there is a path
$$\tau^n Y = A_0 \to A_1 \to \cdots \to A_m = X.$$
Note that such a path is necessarily sectional and $n \leq m$.  By turning the first $n$ arrows one gets a path from $A_n$ to $Y$ and a path from $A_n$ to $X$ which are sectional by minimality of $n$.  Hence let $Z=A_n$.

First, we will prove $\Hom(X,Y[-1])=0$.  Considering the non-split triangle $Y[-1] \to M \to Z \to Y$ we see there is a path from $Y[-1]$ to $Z$.  If $\Hom(X,Y[-1])\not=0$, then there would be a path $Z \to X \to Y[-1]$, contradicting directedness.

Next we will consider $\Hom(X,Y[1])=0$.  We have $\Hom(X,Y[1]) \cong \Hom(Y,\tau X)^{*}$ where the last is shown to be 0 as in the proof of Lemma \ref{lemma:ComponentStandard}.

Finally, we will show $\Hom(X,Y[2])=0$.  If $\Hom(X,Y[2])\not=0$, then $\Hom(Y[1],\tau X)\not=0$.  Since $Z$ has sectional paths to both $X$ and $Y$, we have $\Hom(Z,X)^{*}\cong\Hom(\tau X, \tau^{2} Z[1])\not=0$ and $\Hom(Z,Y)\cong\Hom(Z[1],Y[1])\not=0$.  This gives morphisms $Z[1] \to Y[1] \to \tau X \to \tau^{2} Z[1]$, contradicting directedness.
\end{proof}

Having proved Theorem \ref{theorem:Components}, we now turn our attention to the possible shapes of the Auslander-Reiten components.  First we will discuss a tool we will be developing and using in the next sections.  

\subsection{A word on probing}\label{section:Probing}
In this section, we describe the probing technique.  Following results summarize the results of $\S\ref{subsection:ZDynkin}$, $\S\ref{subsection:ZAi}$, $\S\ref{subsection:ZAii}$, and $\S\ref{subsection:ZDi}$.

As usual, $\AA$ is a connected directed abelian hereditary $k$-linear Ext-finite category satisfying Serre duality, and we write $\CC=\Db\AA$.  We have proven in Theorem \ref{theorem:Components} that the only occurring Auslander-Reiten quivers are of the form $\bZ Q$ where $Q$ is either $A_{\infty}$, $A_{\infty}^{\infty}$, $D_{\infty}$ or a Dynkin quiver.  In Proposition \ref{proposition:ZDynkin} will be proven that if $\CC$ has an Auslander-Reiten component $\bZ Q$ where $Q$ is a Dynkin quiver, then this is the only component of $\CC$.  Since we are interested in the connection between different components, we will exclude such Auslander-Reiten components from this section.  

We will start our discussion with a definition.

Let $\UU$ and $\UU'$ be Auslander-Reiten components.  We will say \emph{$\UU$ maps to $\UU'$} if there is an object $X \in \UU$ and $Y \in \UU'$ such that $\Hom(X,Y) \not= 0$.

It will turn out that the $\ZAi$-components, also called \emph{wings}, are the building blocks of the category $\CC$.  We consider the following map.
\begin{eqnarray*}
\phic:\left\{\mbox{components of $\CC$}\right\} &\to& \left\{ \mbox{sets of wings of $\CC$} \right\} \\
\UU &\mapsto& \{ \mbox{$\WW \mid \WW$ is a wing that maps to $\UU$} \}
\end{eqnarray*}

We now prove some properties of $\phic$.

\begin{proposition}\label{proposition:ProbingComponents}
The map $\phic$ is injective.  Also
\begin{itemize}
\item if $\UU$ is a $\ZAi$-component, then $\phic(\UU) = \{ \UU[-1], \UU\}$,
\item if $\UU$ is a $\ZAii$-component, then $\phic(\UU) = \{ \VV, \WW\}$, with $\VV \not= \WW[z]$ for all $z\in \bZ$,
\item if $\UU$ is a $\ZDi$-component, then $\phic(\UU)$ consists of a single wing.
\end{itemize}
\end{proposition}

\begin{proof}
If $\UU$ is a $\ZAi$-component, a $\ZAii$-component, or a $\ZDi$-component, then the shape of $\phic$ is a direct consequence of Propositions \ref{proposition:LinkmapAi} and \ref{proposition:NoMorphismsBetweenZAi}, Propositions \ref{proposition:LinkmapAii} and \ref{proposition:2UniqueQSinAii}, or Propositions \ref{proposition:LinkmapDi} and \ref{proposition:2UniqueQSinDi}, respectively.  Injectivity of $\phic$ follows from Propositions \ref{proposition:NoMorphismsBetweenZAi}, \ref{proposition:LinkMapAiiInjective} and \ref{proposition:LinkMapDiInjective}.
\end{proof}

We now turn our attention from the components to the objects.  Again, we start with a definition.

We will say that an indecomposable object $S \in \ind \CC$ is a {\em peripheral} object if the middle term of the right Auslander-Reiten triangle is indecomposable.  A peripheral object lying in a wing is a {\em quasi-simple} object.

Quasi-simple objects will be used to, in a certain sense, give coordinates to objects of $\CC$ much like wings can be used as coordinates for components.  We define the function
\begin{eqnarray*}
\phio:\left\{\mbox{indecomposables of $\CC$}\right\} &\to& \left\{ \mbox{sets of quasi-simples of $\CC$} \right\} \\
X &\mapsto& \{ \mbox{$S \mid S$ is a quasi-simple that maps non-zero to $X$} \}
\end{eqnarray*}

\begin{proposition}\label{proposition:ProbingObjects}
Let $X$ be an indecomposable object lying in an Auslander-Reiten component $\UU$.  We have the following properties. 
\begin{enumerate}
\item\label{proposition:ProbingObjects:FromEveryWing} For all $\WW \in \phic(\UU)$, there is an $S \in \phio(X)$ such that $S \in \WW$.
\item\label{proposition:ProbingObjects:Kardinality} The set $\phio(X)$ consists of two elements, except if $X$ is a peripheral object from a $\ZDi$-component, then $\phio(X)$ has only one element.
\item\label{proposition:ProbingObjects:Fiber} The fiber of $\phio(X)$ consists of one element, except when $X$ is a peripheral object in a $\ZDi$-component, then the fiber of $\phio(X)$ consists of two elements.
\item\label{proposition:ProbingObjects:DimensionOne} If $S \in \phio(X)$ then $\dim \Hom(S,X) = 1$.
\item\label{proposition:ProbingObjects:IrreducibleCone} If $S$ is a quasi-simple and $f:S\to X$ a non-zero non-invertible morphism, then the map $g$ in the triangle $\xymatrix@1{S\ar[r]^{f}&X\ar[r]^{g}&C\ar[r]\ar[r]&S[1]}$ is irreducible, except if $X$  is a peripheral object from a $\ZDi$-component.

\end{enumerate}
\end{proposition}

\begin{proof}
\begin{enumerate}
\item
This follows from Propositions \ref{proposition:LinkmapAi}, \ref{proposition:LinkmapAii} and \ref{proposition:LinkmapDi}.
\item
First assume that $X$ is not a peripheral object from a $\ZDi$-component.  Propositions \ref{proposition:LinkmapAi}, \ref{proposition:LinkmapAii} and \ref{proposition:LinkmapDi} yield that there are at least two different quasi-simple objects mapping non-zero to $X$. Propositions \ref{proposition:NoMorphismsBetweenZAi}, \ref{proposition:2UniqueQSinAii} and \ref{proposition:2UniqueQSinDi} imply that these are unique.

If $X$ is a peripheral object from a $\ZDi$-component, then Proposition \ref{proposition:LinkmapDi} yields that there is at least one quasi-simple object mapping non-zero to $X$.  Finally, Proposition \ref{proposition:2UniqueQSinDi} then shows this quasi-simple is unique.
\item
Again, first assume $X$ is not a peripheral object of a $\ZDi$-component.

If $\UU$ is a $\ZAi$-component, then Proposition \ref{proposition:ProbingComponents} yields $\phic(\UU) = \{ \UU[-1], \UU \}$.  We may infer from (\ref{proposition:ProbingObjects:FromEveryWing}) and (\ref{proposition:ProbingObjects:Kardinality}) that $\phio(X) = \{S,T\}$ with $S \in \UU[-1]$ and $T \in \UU$.  Proposition \ref{proposition:LinkmapAi} now yields that the restricted function $\phio\mid_{\ZAi}$ is injective.

If $\UU$ is a $\ZAii$-component, then Proposition \ref{proposition:ProbingComponents} yields $\phic(\UU) = \{ \VV, \WW \}$, with $\VV \not= \WW[z]$ for all $z\in \bZ$.  Now, (\ref{proposition:ProbingObjects:FromEveryWing}) and (\ref{proposition:ProbingObjects:Kardinality}) yield that $\phio(X) = \{S,T\}$ with $S \in \VV$ and $T \in \WW$.  By Proposition \ref{proposition:LinkmapAii} we see that the restricted function $\phio\mid_{\ZAii}$ is injective.

If $\UU$ is a $\ZDi$-component, then Proposition \ref{proposition:ProbingComponents} yields $\phic(\UU) = \{ \VV \}$ and (\ref{proposition:ProbingObjects:FromEveryWing}) and (\ref{proposition:ProbingObjects:Kardinality}) imply that $\phio(X) = \{ S,T \}$ with $S,T \in \VV$.  We may now use Proposition \ref{proposition:LinkmapDi} to see that the restricted function $\phio\mid_{\ZDi}$ is injective.

We may now conclude that the fiber of $\phio(X)$ consists of only one object, $X$.

Now, assume $X$ is a peripheral object of a $\ZDi$-component $\UU$.  We have already shown that there is a unique quasi-simple object, $S$, such that $\Hom(S,X) \not= 0$ and $S \in \WW$ where $\phic(\UU)=\{ \WW \}$.  Since $\phic$ is injective, $\UU$ is the only $\ZDi$-component where $\WW$ maps to.  Proposition \ref{proposition:LinkmapDi} now yields that the fiber of $\phio(X)=S$ consists of exactly two objects, both peripheral objects of $\UU$.

\item[(4 \& 5)] These are immediate consequences of Propositions \ref{proposition:LinkmapAi}, \ref{proposition:LinkmapAii} and \ref{proposition:LinkmapDi}.
\end{enumerate}
\end{proof}

Proposition \ref{proposition:ProbingObjects}(\ref{proposition:ProbingObjects:IrreducibleCone}) will be used in combination with following lemma from \cite{ARS}, adapted to the triangulated case.

\begin{lemma}\label{lemma:irredFact}
Consider the triangle $\xymatrix@1{X\ar[r]^{f}&Y\ar[r]^{g}&Z\ar[r]&X[1]}$ and a morphism $h:Z'\to Z$.  If $f$ is irreducible, then $g$ factors through $h$ or vice versa, thus there exists a $t:Z'\to Y$ such that $h=gt$ or an  $s:Y\to Z'$ such that $g=hs$.
\end{lemma}

\begin{proof}
Consider the morphism of triangles
$$\xymatrix{X\ar[r]^{u}\ar@{=}[d]&C\ar[r]^{v}\ar[d]^{w}&Z'\ar[r]\ar[d]^{h}&X[1]\ar@{=}[d] \\
X\ar[r]_{f}&Y\ar[r]_{g}&Z\ar[r]&X[1]
}$$
Because $f$ is irreducible, we know that $u$ is split mono (and thus $v$ split epi) or $w$ split epi.  In the former case there exists a $t:Z'\to Y$ such that $h=gt$ while in the latter there is a morphism $s:Y\to Z'$ such that $g=hs$.
\end{proof}

\begin{example}
Let $\LL = \bZ$ and consider the category $\AA = \repc(\AL)$ as in \S\ref{section:AL}.  As usual, we write $\CC = \Db \AA$.  Consider a non-zero $f \in \Ext(A_{-1,1},A_{0,3})$ and the triangle $A_{0,3} \to M \to A_{-1,1} \stackrel{f}{\rightarrow} A_{0,3}[1]$.  We will probe $M$ to identify the direct summands of $M$.

One has $\phio(A_{-1,1})=\{ A_{1,1}, A_{-2,-2}[-1] \}$ and $\phio(A_{0,3})=\{ A_{3,3}, A_{-1,-1}[-1] \}$ and may easily verify that the triangle extended with all quasi-simple objects is
$$\xymatrix@C-9pt{A_{1,1}[-1]\ar[d]&A_{3,3}\ar[d]\ar@{=}[r]&A_{3,3}\ar[dr]&&A_{1,1}\ar[dl]&A_{1,1}\ar[d]\ar@{=}[l]&A_{3,3}[1]\ar[d]\\
A_{-1,1}[-1]\ar[r]&A_{0,3}\ar[rr]&&M\ar[rr]&&A_{-1,1}\ar[r]^{f}&A_{0,3}[1]\\
A_{-2,-2}[-2]\ar[u]&A_{-1,-1}[-1]\ar[u]\ar@{=}[r]&A_{-1,-1}[-1]\ar[ur]&&A_{-2,-2}[-1]\ar[ul]&A_{-2,-2}[-1]\ar[u]\ar@{=}[l]&A_{-1,-1}\ar[u]
}$$
We know that $\CC$ does not have any $\ZDi$-components, and thus, by Proposition \ref{proposition:ProbingObjects}(\ref{proposition:ProbingObjects:Kardinality}), that all objects have exactly two quasi-simples mapping non-zero to them.  Since there are exactly four quasi-simples mapping to $M$, we may conclude that $M$ has exactly two direct summands, $M_{1}$ and $M_{2}$.

Using Proposition \ref{proposition:ProbingObjects}(\ref{proposition:ProbingObjects:FromEveryWing}), we may infer that there are two possibilities, either $\phio(M_{1})=\{ A_{1,1}, A_{-2,-2}[-1] \}$ and $\phio(M_{2})=\{ A_{3,3}, A_{-1,-1}[-1] \}$, or $\phio(M_{1})=\{ A_{1,1}, A_{-1,-1}[-1] \}$ and $\phio(M_{2})=\{ A_{3,3}, A_{-2,-2}[-1] \}$.

In the former case, Proposition \ref{proposition:ProbingObjects}(\ref{proposition:ProbingObjects:Fiber}) yields $M_{1} \cong A_{-1,1}$ and $M_{2} \cong A_{0,3}$.  Lemma \ref{lemma:RingelTriangles} then implies that there exists a non-zero morphism from $A_{-1,1}$ to $M_{1}$.  This morphism is necessarily an isomorphism; we conclude that $A_{-1,1} \to M$ is a split monomorphism and hence that $f=0$.  A contradiction.

In the latter case, \ref{proposition:ProbingObjects}(\ref{proposition:ProbingObjects:Fiber}) yields $M_{1} \cong A_{-1,3}$ and $M_{2} \cong A_{0,1}$.

\end{example}

\subsection{A $\bZ Q$-component with $Q$ a Dynkin quiver}\label{subsection:ZDynkin}

We first consider a category $\CC = \Db\AA$ whose Auslander-Reiten quiver has a $\bZ Q$-component where $Q$ is a Dynkin quiver.  Note that the categories $\rep(Q)$ and $\repc(Q)$ are equivalent.

Following proposition shows we may exclude these components from our further discussion of the other components.

\begin{proposition}\label{proposition:ZDynkin}
Let $\AA$ be a connected directed hereditary abelian $k$-linear Ext-finite category satisfying Serre duality.  Assume the Auslander-Reiten quiver of $\CC = \Db\AA$ has a $\bZ Q$-component with $Q$ a Dynkin quiver, then $\CC \cong \Db(\rep (Q^\circ))$.
\end{proposition}

\begin{proof}
Proposition \ref{proposition:SectionPTS} yields that the section $Q$ in the $\bZ Q$-component $\KK$ is a partial tilting set.  Using the exactness of the Serre functor $F$, it is easily seen that the full and exact embedding $\Delta : \Db(\rep(Q^\circ)) \r \CC$ given by Theorem \ref{theorem:PTS} commutes with Serre duality.  Thus those indecomposable objects in the essential image of $\Delta$ are exactly those whose isomorphism class lie in $\KK$.
We claim that for all $X$ in the essential image of $\Delta$ and for all $Y \in \Ob \CC$, we have $\infrad(X,Y) = 0$.  Indeed, if $\infrad(X,Y) \not= 0$ then, for all $n \in \bN$, there would be an $X_{n} \in \KK$ such that $\radn(X,X_{n}) \not= 0$.  Yet, this is not true in $\Db(\rep(Q^{\circ}))$.

Hence, since $\AA$ is connected, $\CC$ can consist of only one component.  We conclude $\CC \cong \Db(\rep(Q^\circ))$.
\end{proof}

We give a further result in this context.

\begin{proposition}\label{proposition:DEquivalentToDynkin}
Let $\AA$ be a hereditary category.  If $\Db\AA \cong \Db(\repc(Q))$ where $Q$ is a Dynkin quiver, then $\AA \cong \rep(Q')$ where the quiver $Q'$ is a tilt of $Q$.
\end{proposition}

\begin{proof}
First note that $\Db(\repc(Q))$ and thus also $\Db\AA$ are directed and $\Ext$-finite.  Since $\AA$ and $\repc(Q)$ are hereditary, we have
\begin{eqnarray*}
\card (\ind \AA) &=& \card \{ \mbox{$T$-orbits in $\ind \Db\AA$} \} \\
								 &=& \card \{ \mbox{$T$-orbits in $\ind \Db(\repc(Q))$} \} \\
								 &=& \card (\ind \repc(Q))
\end{eqnarray*}
Hence $\ind \AA$ is finite.  We will now show this implies that $\AA$ has enough projectives.

Indeed, let $X \in \Ob(\AA)$.  If $X$ is not projective, there exists an object $M = \oplus_{i} M_{i}$  where $M_{i}$ is indecomposable for all $i$, such that there is a non-split epimorphism $M \to X$.  Since $\AA$ is directed and $\ind \AA$ is finite, every sequence of non-split epimorphisms $\cdots \to M \to X$ needs to be finite and we deduce the existence of a projective object $P$ that admits a non-split epimorphism $P \to X$.

Consider the object $P = \oplus_{j} P_{j}$ where $P_{j}$ ranges through all projectives of $\ind \AA$.  We see that $P$ is a generator, and hence $\AA \cong \mod(A)$ with $A = \End(P)$.

Since $A$ is hereditary and of finite representation type, $A$ needs to be Morita equivalent to the path algebra of a Dynkin quiver $Q'$.

Finally, note that the Auslander-Reiten quiver of $\Db \AA$ and $\Db \repc(Q)$, are equal to $\bZ Q'$ and $\bZ Q$, respectively, hence $Q'$ is a tilt of $Q$.
\end{proof}
\subsection{A $\ZAi$-component}\label{subsection:ZAi}

\begin{figure}
	\centering
		\includegraphics{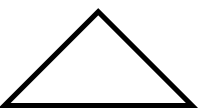}
	\caption{A $\ZAi$-component}
	\label{fig:ZAi}
\end{figure}

A $\ZAi$-component will often be called a wing and when we encounter such a component in a derived category, we will represent it by a triangle as shown in Figure \ref{fig:ZAi}.

\begin{proposition}\label{proposition:LinkmapAi}
Let $Q$ be the quiver
$$\xymatrix{
0\ar[rd]&&2\ar[ld]\ar[rd]&&4\ar[ld]\ar@{.}[rd] \\
&1&&3&&
}$$
and let $\KK$ be a $\ZAi$-component of $\CC$, then the smallest full and exact subcategory of $\CC$ containing $\KK$ is equivalent to $\Db \rep Q$, and the embedding $\Delta : \Db \rep Q \to \CC$ commutes with Serre duality.  Hence $\Delta$ maps Auslander-Reiten components to Auslander-Reiten components.
\end{proposition}

\begin{proof}
Consider within $\KK$ the quiver $Q^{\circ}$ as in Figure \ref{Fig:QinZAi}.  We will denote the indecomposable corresponding to the vertex $i$ of $Q^{\circ}$ by $P_{i}$.

\begin{figure}
	\centering
$$\xymatrix@C=20pt@R=20pt{&{\vdots}\ar@{..>}[dr]&&{\vdots}\ar@{..>}[dr]&&\vdots\ar@{..>}[dr]&&{\vdots}\ar@{..>}[dr]\\
{\cdots}\ar@{..>}[dr]\ar@{..>}[ur]&&{\bullet}\ar@{..>}[dr]\ar@{..>}[ur]&&P_{3}\ar[dr]\ar[ur]&&{\bullet}\ar@{..>}[dr]\ar@{..>}[ur]&&{\cdots}\\
&{\bullet}\ar@{..>}[dr]\ar@{..>}[ur]&&{\bullet}\ar@{..>}[ur]\ar@{..>}[dr]&&P_{2}\ar@{..>}[dr]\ar@{..>}[ur]&&{\bullet}\ar@{..>}[ur]\ar@{..>}[dr]\\
{\cdots}\ar@{..>}[dr]\ar@{..>}[ur]&&{\bullet}\ar@{..>}[dr]\ar@{..>}[ur]&&P_{1}\ar[dr]\ar[ur]&&{\bullet}\ar@{..>}[dr]\ar@{..>}[ur]&&{\cdots}\\
&{\bullet}\ar@{..>}[ur]&&{\bullet}\ar@{..>}[ur]&&P_{0}\ar@{..>}[ur]&&{\bullet}\ar@{..>}[ur]
}$$
	\caption{The quiver $Q^{\circ}$ in a $\ZAi$-component}
	\label{Fig:QinZAi}
\end{figure}
Invoking Proposition \ref{proposition:SectionPTS} and Theorem \ref{theorem:PTS}, we may consider a full and exact embedding $\Delta : \Db \rep Q \to \CC$ which we claim to commute with the Serre functor.

Considering the exactness of $\Delta$, and the connection between the Auslander-Reiten translation $\tau$ and the Serre functor $F$, it is easy to see $\Delta F P_{i}\cong F \Delta P_{i}$ for all $i\in \bN$.  Since the Serre functor is exact and commutes with $\Delta$ on generators of $\Db \rep Q$, it will commute with $\Delta$.

We still need to check whether $\Delta$ maps Auslander-Reiten components to Auslander-Reiten components essentially surjective, i.e.\ if an indecomposable of a component is in the essential image of $\Delta$, then so is every indecomposable of that component.  To this end, consider an indecomposable object $C$ in $\CC$ in the essential image of $\Delta$ such that there is an irreducible $D \to C$ where $D$ is not in the essential image of $\Delta$ (the dual case where there is an irreducible $C \to D$ is completely analogue).  If $C$ is in the essential image of $\Delta$, then so is $\tau C$ since $\Delta$ commutes with $F$.  We may consider the Auslander-Reiten triangles
$$\xymatrix{
\tau M \ar[r]\ar@{.>}[d]^{\Delta}& N' \ar[r]\ar@{.>}[d]^{\Delta}& M \ar[r]\ar@{.>}[d]^{\Delta}&\tau M[1]\ar@{.>}[d]^{\Delta} \\
\tau C \ar[r]&N \ar[r]& C \ar[r]&A[1]}$$
Since $\Delta$ is full, faithful and exact, we know that $\End(N) \cong \End(N')$, hence $N$ and $N'$ consists of the same number of indecomposable summands and there must be an indecomposable direct summand $D'$ of $N'$ such that $\Delta(D') \cong D$.
\end{proof}

\begin{remark}
The category $\rep Q$ occurring in the proof has been described in \cite{ReVdB02}. We may sketch the bounded derived category $\Db \rep Q$ as shown in Figure \ref{fig:DbAi}
\begin{figure}
	\centering
		\includegraphics{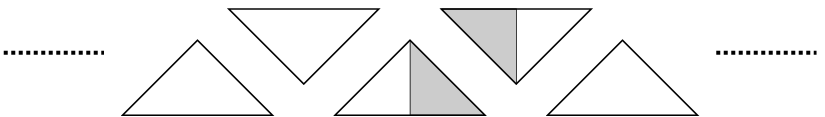}
	\caption{A sketch of $\Db\rep Q$ occurring in Proposition \ref{proposition:LinkmapAi}}
	\label{fig:DbAi}
\end{figure}
where we have marked the abelian subcategory $\rep Q$ with gray.  Note that $\Db \rep Q \cong \Db \repc \bZ$.
\end{remark}

\begin{proposition}\label{proposition:NoMorphismsBetweenZAi}
Let $\VV$ and $\WW$ be wings.  If $\VV$ maps to $\WW$, then $\WW = \VV$ or $\WW = \VV[1]$.
\end{proposition}

\begin{proof}
It is clear that $\VV$ maps to $\WW$ if $\WW = \VV$ or $\WW = \VV[1]$.

To prove that $\VV$ does not map to $\WW$ otherwise, we start by fixing a notation.  Let $Q$ be the quiver
$$0 \to 1 \to 2 \to 3\to \cdots$$
Since $\VV \cong \bZ Q \cong \WW$ as stable translation quivers, we may label the vertices of $\VV$ by $\VV_{m,n}$ and the vertices of $\WW$ by $\WW_{m,n}$ for $m \in \bZ$ and $n \in \bN$ as is illustrated in Figure \ref{Fig:LabelingVV}.

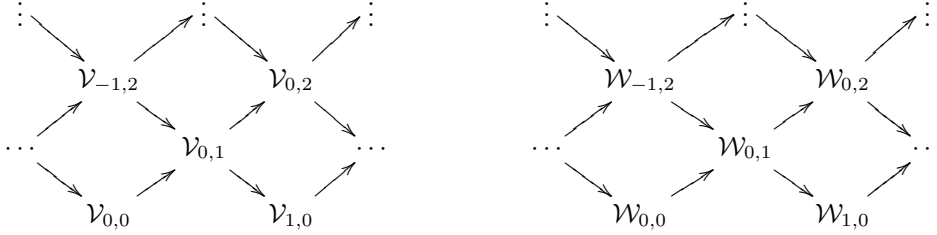
\begin{figure}
\exa
$\xymatrix@C=10pt@R=10pt{{\vdots}\ar[dr]&&{\vdots}\ar[dr]&&\vdots\\
&\VV_{-1,2}\ar[dr]\ar[ur]&&\VV_{0,2}\ar[ur]\ar[dr]\\
\cdots\ar[dr]\ar[ur]&&\VV_{0,1}\ar[dr]\ar[ur]&&\cdots\\
&\VV_{0,0}\ar[ur]&&\VV_{1,0}\ar[ur]
}$
\exb
$\xymatrix@C=10pt@R=10pt{{\vdots}\ar[dr]&&{\vdots}\ar[dr]&&\vdots\\
&\WW_{-1,2}\ar[dr]\ar[ur]&&\WW_{0,2}\ar[ur]\ar[dr]\\
\cdots\ar[dr]\ar[ur]&&\WW_{0,1}\ar[dr]\ar[ur]&&\cdots\\
&\WW_{0,0}\ar[ur]&&\WW_{1,0}\ar[ur]
}$
\exc
	\caption{Labeling the vertices in the components $\VV$ and $\WW$ from Proposition \ref{proposition:NoMorphismsBetweenZAi}}
	\label{Fig:LabelingVV}
\end{figure}

One sees easily that
$$\dim\Hom(\VV_{m,n},\WW_{i,1}) = \dim\Hom(\VV_{m,n},\WW_{i,0})+\dim\Hom(\VV_{m,n},\WW_{i+1,0})$$
and, by induction, that
$$\dim\Hom(\VV_{m,n},\WW_{i,j}) = \sum_{k=0}^{j}\dim\Hom(\VV_{m,n},\WW_{i+k,0})$$
Dually, one has
$$\dim\Hom(\VV_{m,n},\WW_{i,j}) = \sum_{l=0}^{n}\dim\Hom(\VV_{m+l,0},\WW_{i,j}).$$
Assume $\Hom(\VV_{m,n},\WW_{i,j})\not=0$, then there are $p$ and $q$ in $\bN$ with $m\leq p \leq m+n$ and $i\leq q \leq i+j$ such that $\Hom(\VV_{p,0},\WW_{q,0})\not=0$, and hence $\Hom(\tau^{k} \VV_{p,0},\tau^{k} \WW_{q,0}) \cong \Hom(\VV_{p-k,0},\WW_{q-k,0})\not=0$ for all $k\in \bZ$.  Finally, we obtain
\begin{eqnarray*}
\dim\Hom(\VV_{p,6},\WW_{q,6})&=&\sum_{k=0}^{6}\dim\Hom(\VV_{p,6},\WW_{q+k,0})\\
&=&\sum_{k,l=0}^{6}\dim\Hom(\VV_{p+l,0},\WW_{q+k,0})\\
&\geq&\sum_{k=0}^{6}\dim\Hom(\VV_{p+k,0},\WW_{q+k,0})\\
&\geq&7
\end{eqnarray*}
contradicting Theorem \ref{theorem:BoundedHomDim}.
\end{proof}
\subsection{A $\ZAii$-component}\label{subsection:ZAii}

\begin{figure}
	\centering
		\includegraphics{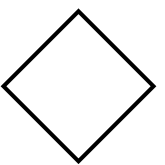}
	\caption{A $\ZAii$-component}
	\label{fig:Diamond}
\end{figure}

In this part, we will discuss the $\ZAii$-component.  A $\ZAii$-component in a derived category will be drawn as in Figure \ref{fig:Diamond}.

\begin{proposition}\label{proposition:LinkmapAii}
Let $Q$ be the quiver 
$$\xymatrix{
\ar@{.}[rd]&&-2\ar[ld]\ar[rd]&&0\ar[ld]\ar[rd]&&2\ar[ld]\ar[rd]&&4\ar[ld]\ar@{.}[rd] \\
&-3&&-1&&1&&3&&
}$$
and let $\KK$ be a $\ZAii$-component of $\CC$, then the smallest full and exact subcategory of $\CC$ containing $\KK$ is equivalent to $\Db \rep Q$, and the embedding $\Delta : \Db \rep Q \to \CC$ commutes with Serre duality.  Hence $\Delta$ maps Auslander-Reiten components to Auslander-Reiten components.
\end{proposition}

\begin{proof}
The construction of the functor $\Delta$ is similar to the construction in the proof of Proposition \ref{proposition:LinkmapAi}; one now finds the quiver $Q^\circ$ within the $\ZAii$-component $\KK$. 
\end{proof}

\begin{remark}
The category $\rep Q$ occurring in the proof has been discussed in \cite{ReVdB02}.  We may sketch the derived category $\Db \rep Q$ as in Figure \ref{fig:Aii} where we have filled the abelian subcategory $\rep Q$ with gray.
\begin{figure}
	\centering
		\includegraphics{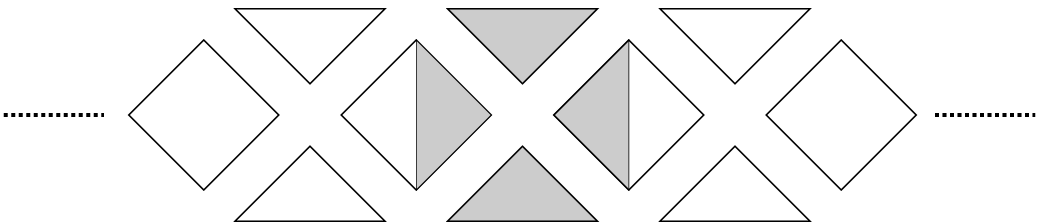}
	\caption{A sketch of $\Db \rep Q$ occuring in Proposition \ref{proposition:LinkmapAii}}
	\label{fig:Aii}
\end{figure}
Note that the category $\repc (\{0,1\} \stackrel{\rightarrow}{\times} \bZ)$ is derived equivalent to $\rep Q$.
\end{remark}

Since the embedding $\Delta : \Db \rep Q \to \CC$ maps Auslander-Reiten components to Auslander-Reiten components, we may define the map 
$$\phi_{\ZAii}:\left\{\mbox{$\ZAii$-components of $\CC$}\right\} \to\left\{ \mbox{pairs of $\ZAi$-components of $\CC$} \right\}$$
by mapping a $\ZAii$-components $\KK$ to the pair of wings within the essential image of $\Delta$ that map non-zero to $\KK$.

\begin{remark}
Note that the map $\phi_{\ZAii}$ does not depend on the choice of the partial tilting set in the proof of Proposition \ref{proposition:LinkmapAii}.
\end{remark}

\begin{remark}
It will follow from the following Propositions that $\phi_{\ZAii} = \left. \phic \right|_{\ZAii}$.
\end{remark}

\begin{proposition}\label{proposition:2UniqueQSinAii}
If a wing $\WW$ maps to a $\ZAii$-component $\KK$, then $\WW \in \phi_{\ZAii}(\KK)$.
\end{proposition}

\begin{proof}
We start by fixing an $X \in \KK$.  An argument analogous to the proof of Proposition \ref{proposition:NoMorphismsBetweenZAi} shows that there are only finitely many quasi-simple objects of $\WW$ that map non-zero to $X$.  We may choose a quasi-simple $S$ from $\WW$ such that $\Hom(\tau^{-1}S,X)=0$ and $\Hom(S,X)\not=0$.  Now, consider the triangle $\tri{S}{X}{C}$.  Applying the functor $\Hom(X,-)$, we may conclude easily that $\dim\Hom(X,C)=1$.

Next, consider the Auslander-Reiten triangle $\tri{X}{Y\oplus Y'}{\tau^{-1}X}$.  Since the morphism $X\to C$ factors through $Y \oplus Y'$, we may assume there exists a morphism $X \to Y$ such that the composition $X\to Y\to C$ is non-zero.  This gives rise to the following morphism of triangles
$$\xymatrix{
S\ar[r]&X\ar[r]&C\ar[r]&S[1]\\
T\ar[r]\ar[u]&X\ar[r]\ar@{=}[u]&Y\ar[u]\ar[r]&T[1]\ar[u]
}$$
Since $X \to Y$ is an irreducible morphism between indecomposable objects, Proposition \ref{proposition:LinkmapAii} yields that $T$ is a quasi-simple object from a wing $\VV \in \phi_{\ZAii}(\KK)$.  The induced morphism $T \to S$ is easily seen to be non-zero.  Proposition \ref{proposition:NoMorphismsBetweenZAi} yields $\WW = \VV$ or $\WW = \VV[1]$.  By Proposition \ref{proposition:LinkmapAii} we may exclude the latter.  We conclude $\WW \in \phi_{\ZAii}(\KK)$.
\end{proof}





\begin{proposition}\label{proposition:LinkMapAiiInjective}
The map $$\phi_{\ZAii}:\left\{\mbox{$\ZAii$-components of $\CC$}\right\} \to\left\{ \mbox{pairs of $\ZAi$-components of $\CC$} \right\}$$ is injective.
\end{proposition}

\begin{proof}
Let $\KK$ en $\KK'$ be $\ZAii$-components such that $\phi_{\ZAii}(\KK) = \{\VV,\WW\} = \phi_{\ZAii}(\KK')$.  Fix a quasi-simple object $S$ from the component $\VV$ and a quasi-simple object $T$ from the component $\WW$.  Proposition \ref{proposition:LinkmapAii} yields a unique indecomposable object $X$ from $\KK$ and a unique indecomposable object $X'$ from $\KK'$ such that $\Hom(S,X), \Hom(T,X), \Hom(S,X')$ and $\Hom(T,X')$ are all non-zero.  Furthermore, it follows from Proposition \ref{proposition:LinkmapAii} that all these Hom-spaces are 1-dimensional and from Proposition \ref{proposition:2UniqueQSinAii} that these are the only quasi-simples mapping to $X$ or $X'$.  We wish to prove that $X \cong X'$.

Since by Proposition \ref{proposition:LinkmapAii} the map $X\to Y_{1}$ occurring in the triangle $\tri{S}{X}{Y_{1}}$ is irreducible, we may use Lemma \ref{lemma:irredFact} to see there is a morphism $X\to X'$ or $X'\to X$.  Thus without loss of generality, we may conclude there is a commutative diagram
$$\xymatrix{S\ar[r]&X'\\
S\ar[r]\ar@{=}[u]&X\ar[u]
}$$

Analogously, considering the triangle $\tri{T}{X}{Y_{2}}$ and directedness gives a morphism $X\to X'$ and we obtain the commuting diagram
$$\xymatrix{T\ar[r]&X'\\
T\ar[r]\ar@{=}[u]&X\ar[u]
}$$

It easily follows there is a morphism $f:X\to X'$ such that both compositions $S\to X\to X'$ and $T\to X \to X'$ are non-zero.

In order to prove that $f$ is an isomorphism, we will use quasi-simples to \emph{probe} the object $M = \cone(f:X \to X')$, i.e.\ we will look which quasi-simple objects map to $M$.  The triangle built on $f$, extended with all the quasi-simple objects mapping to each of its objects looks like
$$\xymatrix{S\ar[d]\ar@{=}[r]&S\ar[d]&&S[1]\ar[d]\ar@{=}[r]&S[1]\ar[d]\\
X\ar[r]^{f}&X'\ar[r]^{g}&M\ar[r]^{h}&X[1]\ar[r]^{f[1]}&X'[1]\\
T\ar[u]\ar@{=}[r]&T\ar[u]&&T[1]\ar[u]\ar@{=}[r]&T[1]\ar[u]
}$$

We will show that no quasi-simple object $U$ may map to $M$.  Seeking a contradiction, assume that $f' : U \to M$ is such a non-zero morphism.  We will first consider the case where the composition $U \to M \stackrel{h}{\rightarrow} X[1]$ is zero.  In this case the map $U \to M$ factors though $g:X' \to M$, and $U$ would be isomorphic to either $S$ or $T$.  But since then $\dim \Hom(U,X') = 1$, we may further conclude that $U \to X'$ factors through $f:X \to X'$ so that the composition $U \to X \stackrel{f}{\rightarrow} X' \stackrel{g}{\rightarrow} M$ is non-zero, a contradiction.

Analogous, if the map $U \to M \stackrel{h}{\rightarrow} X[1]$ is non-zero, then $U$ would map non-zero to $X[1]$ and would hence be isomorphic to either $S[1]$ or $T[1]$.  Again, since then $\dim \Hom(U,X[1]) = 1$, we may conclude that the composition $U \to M \stackrel{h}{\rightarrow} X[1] \stackrel{f[1]}{\rightarrow} X'[1]$ is non-zero, a contradiction.

Using Propositions \ref{proposition:LinkmapAi} and \ref{proposition:LinkmapAii}, and already using Proposition \ref{proposition:LinkmapDi} from the next section, we see that every non-zero object of $\CC$ has at least one quasi-simple object mapping to it.  The cone $M$ thus has to be the zero object, establishing the fact that $X$ and $X'$ are isomorphic, and thus that $\KK = \KK'$.
\end{proof}
\subsection{A $\ZDi$-component}\label{subsection:ZDi}

\begin{figure}
	\centering
		\includegraphics{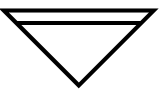}
	\caption{A $\ZDi$-component}
	\label{fig:dtriangle}
\end{figure}

In this part, we will discuss the $\ZDi$-component within the directed category.  Mostly, the proofs are analogous to the case of a $\ZAii$-component.  We will draw a $\ZDi$-component of a derived category as a triangle with a doubled base as in Figure \ref{fig:dtriangle}

\begin{proposition}\label{proposition:LinkmapDi}
Let $Q$ be the quiver 
$$\xymatrix{
&2\ar[rd]\ar[d]\ar[ld]&&4\ar[ld]\ar[rd]&&6\ar[ld]\ar[rd]&&8\ar[ld]\ar@{.}[rd] \\
0&1&3&&5&&7&&
}$$
and let $\KK$ be a $\ZDi$-component of $\CC$, then the smallest full and exact subcategory of $\CC$ containing $\KK$ is equivalent to $\Db\rep Q$, and the embedding $\Delta : \Db \rep Q \to \CC$ commutes with Serre duality.  Hence $\Delta$ maps Auslander-Reiten components to Auslander-Reiten components.
\end{proposition}

\begin{proof}
The proof is analogue to the proof of Proposition \ref{proposition:LinkmapAi}.  One finds the quiver $Q^\circ$ within the $\ZDi$-component $\KK$.
\end{proof}

\begin{remark}
The category $\rep Q$ occurring in Proposition \ref{proposition:LinkmapDi} has been discussed in \cite{ReVdB02}.  We may sketch the derived category $\Db\rep Q$ as in Figure \ref{fig:Di} where we have filled the abelian subcategory $\rep Q$ with gray.
\begin{figure}
	\centering
		\includegraphics{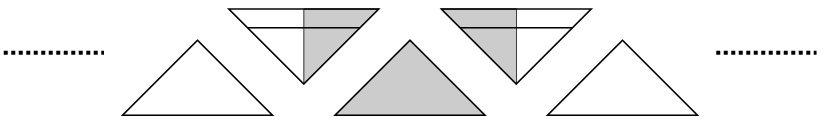}
	\caption{A sketch of $\Db \rep Q$ occuring in Proposition \ref{proposition:LinkmapDi}}
	\label{fig:Di}
\end{figure}
Note that the category $\repc D_\bZ$ is derived equivalent to the category $\rep Q$.
\end{remark}

As in discussion of the $\ZAii$-component, the embedding $\Delta : \Db \rep Q \to \CC$ maps Auslander-Reiten components to Auslander-Reiten components.  Hence, we may define the map 
$$\phi_{\ZDi}:\left\{\mbox{$\ZDi$-components of $\CC$}\right\} \to\left\{ \mbox{singletons of $\ZAi$-components of $\CC$} \right\}$$
by mapping a $\ZDi$-components $\KK$ to the set of wings within the essential image of $\Delta$ that map non-zero to $\KK$.

\begin{remark}
Note that the map $\phi_{\ZDi}$ does not depend on the choice of the partial tilting set in the proof of Proposition \ref{proposition:LinkmapDi}.
\end{remark}

\begin{remark}
It will follow from the following Propositions that $\phi_{\ZDi} = \left. \phic \right|_{\ZDi}$.
\end{remark}

The proofs of the following propositions are analogous to the proofs of the corresponding properties in our discussion of the $\ZAii$-component.

\begin{proposition}\label{proposition:2UniqueQSinDi}
If there exists a non-zero morphism from a $\ZAi$-component $\WW$ to a $\ZDi$-component $\KK$, then $\phi_{\ZDi}(\KK) = \{ \WW \}$.
\end{proposition}

\begin{proposition}\label{proposition:LinkMapDiInjective}
The map $$\phi_{\ZDi}:\left\{\mbox{$\ZDi$-components of $\CC$}\right\} \to\left\{ \mbox{singletons of $\ZAi$-components of $\CC$} \right\}$$ is injective.
\end{proposition}
\section{Classification}

Let $\AA$ be a connected directed hereditary abelian $k$-linear Ext-finite category satisfying Serre duality, and write $\CC=\Db\AA$.  In this section we will prove our main theorem.  We will proceed as follows.  First, we will define a partial tilting set $\St$.  Then Theorem \ref{theorem:PTS} will give a full and exact embedding $\Db(\mod{\St^{\circ}}) \to \Db\AA$.  We will proceed to prove that this embedding is an equivalence of triangulated categories.  Finally, the classification will follow from the shape of the poset $\St$.  We start by defining $\St$.

Choose a quasi-simple object $S$ in a wing $\WW$.  We will consider two cases.  First, assume that $S$ does not map to two peripheral objects $Q^1_{S}$ and $Q_{S}^2$ of a $\ZDi$-component or, equivalently, there is no $\ZDi$-component $\KK$ such that $\WW \in \phi_{\ZDi}(\KK)$.  In this case just let $\St$ be the set of indecomposable objects $X$ such that there exists a map from $S$ to $X$, thus
$$\St = \left\{X \in \ind \CC \mid \Hom(S,X)\not=0\right\}.$$

Secondly, assume $S$ does map to two peripheral objects $Q^1_{S}$ and $Q_{S}^2$ of a $\ZDi$-component or, equivalently, there is a $\ZDi$-component $\KK$ such that $\WW \in \phi_{\ZDi}(\KK)$.  Then, let $\St$ be the set of indecomposable objects $X$ such that there exists a map from $S$ to $X$ and a map from $X$ to $Q_{S}^1$ or $Q_{S}^2$, thus
$$\St = \{X \in \ind \CC |\Hom(S,X)\not=0 \mbox{ and }\Hom(X,Q_{S}^1 \oplus Q_{S}^2)\not=0\}.$$

This defines the full subcategory $\St$ of $\CC$.  We will define a poset structure by 
$$X \leq Y \Leftrightarrow \Hom(X,Y) \not= 0$$
and by identifying isomorphic objects.  Before proving in Lemma \ref{lemma:StPoset} that this does indeed define a poset structure, we will give two examples.

Therefore, we will fix following notation.  Let $\PP_{1}$ and $\PP_{2}$ be posets.  The poset $\PP_{1} \cdot \PP_{2}$ has $\PP_{1} \stackrel{\cdot}{\cup} \PP_{2}$ as underlying set and
$$X \leq Y \Leftrightarrow \left\{ \begin{array}{l} \mbox{$X \leq Y$ in $\PP_{1}$,  or  }\\ \mbox{$X \leq Y$ in $\PP_{2}$,  or  }\\ \mbox{$X \in \PP_{1}$ and $Y \in \PP_{2}$.} \end{array}\right.$$

\begin{example}\label{example:StAL}
Let $\LL$ be the poset $\{0,1\}\stackrel{\rightarrow}{\times}\bZ$.  The abelian category $\repc(\AL)$ consists of two $\ZAi$-components, one containing the indecomposables of the form $A_{(0,i),(0,j)}$ and one containing the indecomposable objects of the form $A_{(1,i),(1,j)}$, and a $\ZAii$-component wherein all the indecomposable objects $A_{(0,i),(1,j)}$ lie, for all $i,j\in\bZ$.  Now, let $S=A_{(1,0),(1,0)}$.  We may then describe the set $\St$ as
$$\St = \{A_{(1,-n),(1,0)}|n\in\bN\} \cdot \{A_{(0,z),(1,0)}|z\in\bZ\} \cdot \{A_{(1,1),(1,n)}[1]|\mbox{$n\in\bN$ and $n\geq 1$}\}.$$
We may draw $\St$ within $\Db(\repc(\AL))$ as in Figure \ref{fig:SAii} where, as usual, the abelian category $\repc(\AL)$ has been filled with gray.
\begin{figure}
	\centering
		\psfrag{S}[][]{$S$}
		\psfrag{ts}[][]{$\tau S[1]$}
		\includegraphics{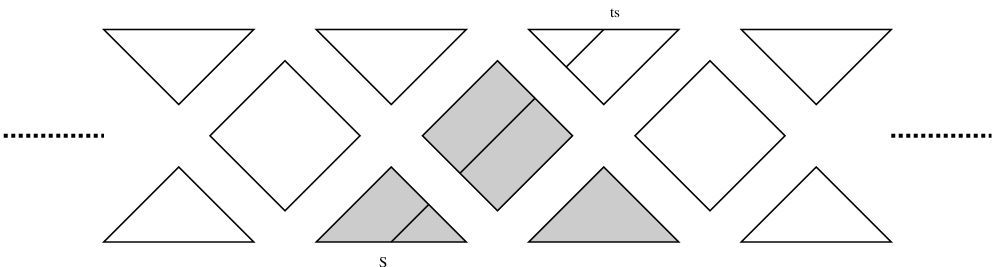}
	\caption{The set $\St$ in $\repc(\AL)$}
	\label{fig:SAii}
\end{figure}
\end{example}

\begin{example}\label{example:StDL}
In this example, we will consider $\DL$ where $\LL = \bZ$.  The abelian category $\repc(\DL)$ consists of a $\ZAi$-component containing the indecomposables of the form $A_{i,j}$ and a $\ZDi$-component containing the indecomposables of the form $B_{i,j}$, $A_{i}^{1}$, and $A_{i}^{2}$.  If $S=A_{0,0}$ then we may describe the set $\St$ as
$$\St =\{A_{-n,0}|n\in\bN\}\cdot\{B_{0,n+1}|n\in\bN\}\cdot\{B_{0}^{1},B_{0}^{2}\}.$$
Graphically, we may represent $\St$ within $\Db(\repc(\LL))$ as in Figure \ref{fig:SDi}.
\begin{figure}
	\centering
		\psfrag{S}[][]{$S$}
		\includegraphics{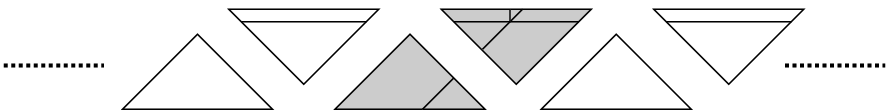}
	\caption{The set $\St$ in $\repc(\DL)$}
	\label{fig:SDi}
\end{figure}
\end{example}

Note that in Examples \ref{example:StAL} and \ref{example:StDL} the set $\St$ is of the form $\AL^{\circ}$ or $\DL^{\circ}$ where $\LL$ is a \emph{bounded} locally discrete linearly ordered set, i.e.\ a locally discrete linearly ordered set with both a minimal and a maximal element.  It is easily seen that for such linearly ordered posets we have $\LL \cong (\bN) \cdot (\TXZ) \cdot (-\bN)$ for a certain linearly ordered set $\TT$.  

Following lemma will classify all possible posets that may occur as $\St$.  For $X,Z \in \St$ we will write
$$[X,Z] = \{ Y \in \St \mid X \leq Y \leq Z \}.$$

\begin{lemma}\label{lemma:StPoset}
The set $\St$ is a poset of the form $\AL^{\circ}$ or $\DL^{\circ}$ where $\LL$ is a bounded locally discrete linearly ordered set.
\end{lemma}

\begin{proof}
We start with the case where $S$ does not map to peripheral objects of a $\ZDi$-component and wish to prove that $\St \cong \AL^{\circ}$ where $\LL$ is a bounded locally discrete linearly ordered set.

The relation $\leq$ in $\St$ is obviously reflexive. The fact that it is antisymmetric follows from directedness.  In order to prove transitivity and linearly ordered, it suffices to prove that $\Hom(X,Y)\not=0$ or $\Hom(Y,X)\not=0$.  For all $X,Y \in \St$, we may consider the commutative diagram
$$\xymatrix{Y\\
S\ar[u]\ar[r]^{f}&X\ar[r]^{g}&C\ar[r]&S[1]}$$
where the bottom line is a triangle and $g:X \to C$ is irreducible as is shown in Proposition \ref{proposition:ProbingObjects}(\ref{proposition:ProbingObjects:IrreducibleCone}).  Lemma \ref{lemma:irredFact} now yields that $\Hom(X,Y)\not=0$ or $\Hom(Y,X)\not=0$, thus $\St$ is a linearly ordered poset.  Even more so, if $X \leq Y$ we may assume that the composition $S \r X \r Y$ is non-zero.

Using Serre duality, it is easily seen that $\St = [S,\tau S[1]]$, thus $\St$ is bounded.

To prove that $\St$ is locally discrete, assume $Z \in \St$ is a non-minimal element. We need to prove there exists a finite set $A \subseteq \St$ such that, for all $X \in \St$ with $X < Z$, there is a $Y \in A$ with $X \leq Y < Z$.

Therefore, consider the Auslander-Reiten triangle $\tri{\tau Z}{M}{Z}$.  Write $M = \oplus_{i} M_{i}$ where $M_{i}$ is indecomposable.  We have already proven that there exists a non-zero morphism from $X$ to $Z$ such that the composition $S \to X \to Z$ is non-zero.  Since $X \to Z$ factors through $M$, it is clear that there exists an $M_{i} \in \ind M$ such that $M_{i} \in \St$ and $X \leq M_{i} < Z$, thus we have shown $A = \St \cap \ind M$.

The case where $Z$ is a non-maximal element is analogous.  Hence $\St$ is locally discrete.



We now turn our attention to the case where $S$ does map to two peripheral objects of a $\ZDi$-component and wish to prove that $\St \cong \DL^{\circ}$ where $\LL$ is a bounded locally discrete linearly ordered set.

Note that the left Auslander-Reiten triangles built on $Q_{S}^1$ and $Q_{S}^2$ both have the same indecomposable middle term $N$.  It is straightforward to check that $\St \cong [S,N]\cdot\{Q_{S}^1, Q_{S}^2\}$.  Analogous to the first part of the proof, one shows that $[S,N]$ is a bounded locally discrete linearly ordered set, thus $\St \cong \DL^{\circ}$ where $\LL = [S,N]^{\circ}$.
\end{proof}

We will now prove most technical results needed to prove our main result, Theorem \ref{theorem:MainTheorem}.  In particular, we will construct an embedding $\Db(\rep(\St^{\circ})) \to \CC$ and prove that it is an equivalence of categories.

\begin{lemma}\label{lemma:StPTS}
The set $\St$ is a partial tilting set.
\end{lemma}

\begin{proof}
We need to prove that $\Hom(X,Y[z])=0$ for all $z \in \bZz$, and all $X,Y\in \St$.  If $X$ and $Y$ are the peripheral objects $Q_{S}^1$ and $Q_{S}^2$ then the assertion follows easily.  Indeed, both are then contained within the same $\ZDi$-component, and all the maps within such component are known from Proposition \ref{proposition:LinkmapDi}.

Since we may now assume that either $X$ or $Y$ is not isomorphic to $Q_{S}^1$ or to $Q_{S}^2$, Lemma \ref{lemma:StPoset} yields that we may assume $\Hom(X,Y) \not= 0$ and hence, by Lemma \ref{lemma:Directed}, we need only to prove that $\Hom(Y,X[1])=0$.

Let $f \in \Hom(Y,X[1])$ and consider the triangle
\begin{equation}\label{eqn:ClassificationTriangle}
\xymatrix@1{X[-1]\ar[r]&Y\ar[r]&M\ar[r]&X\ar[r]^{f}&Y[1]}
\end{equation}
We will now probe $M$ to prove $M \cong X \oplus Y$ and hence $f=0$.

First, consider the case where $X = Q_{S}^1$ and $Y \not = Q_{S}^2$.

Since $Y\in \St$, we have $\Hom(S,Y)\not=0$, and consequently $\Hom(S,Y[1])=0$.  We conclude that the composition $S \to X \stackrel{f}{\rightarrow} Y[1]$ needs to be zero.

Let $S$ and $T'$ be the quasi-simples mapping to $Y$, thus $\phio(Y)=\{S,T'\}$.  Since $\phio(X[-1])=\{S[-1]\}$, we have that $\Hom(S,X[-1])=0$ and $\Hom(T',X[-1])=0$.

If $S'$ is a quasi-simple such that $\Hom(S',M) \not=0$, then either $\Hom(S',X) \not=0$ or $\Hom(S',Y) \not=0$, hence $S' \cong S$ or $S' \cong T'$.

Thus triangle (\ref{eqn:ClassificationTriangle}) enriched with all the quasi-simples mapping to each of its components is

$$\xymatrix{S[-1]\ar[d]&S\ar[d]\ar@{=}[r]&S\ar[dr]&&S\ar[dl]\ar@{=}[r]&S\ar[d]&S[1]\ar[d]\\
X[-1]\ar[r]&Y\ar[rr]&&M\ar[rr]&&X\ar[r]^{f}&Y[1]\\
&T'\ar[u]\ar@{=}[rr]&&T'\ar[u]&&&T'[1]\ar[u]
}$$

Using Proposition \ref{proposition:ProbingObjects} it is easy to see that either $Q_{S}^1$ or $Q_{S}^2$ have to be direct summands of $M$.  However, since all non-zero morphisms in $\Hom(Q^1_{S},Q^1_{S})$ are isomorphisms and $\Hom(Q_{S}^2,Q_{S}^1)=0$, Lemma \ref{lemma:RingelTriangles} implies $f=0$.

We may now assume that neither $X$ nor $Y$ are the peripheral objects $Q_{S}^1$ or $Q_{S}^2$.  Proposition \ref{proposition:ProbingObjects}(\ref{proposition:ProbingObjects:Kardinality}) yields there are two quasi-simple objects, $S$ and $T$, mapping to $X$ and two quasi-simple objects, $S$ and $T'$, mapping to $Y$.  We will have to consider two cases, namely one where the composition $\xymatrix@1{T\ar[r]&X\ar[r]^{f}&Y[1]}$ is zero and one where it is non-zero.  We start with the former.

We will probe $M$ to proof $M \cong X \oplus Y$, so that $f=0$.  Therefore, we wish to find all quasi-simples that admit a non-zero map to $M$.

Since $\Hom(S,Y) \not= 0$, we know by directedness that $\Hom(S,Y[1])=0$ and hence that the composition $\xymatrix@1{S\ar[r]&X\ar[r]^{f}&Y[1]}$ is zero.

If $S'$ is a quasi-simple such that $\Hom(S',M) \not=0$, then either $\Hom(S',X) \not=0$ or $\Hom(S',Y) \not=0$, hence $S' \cong S$, $S' \cong T$, or $S' \cong T'$.

Thus triangle (\ref{eqn:ClassificationTriangle}) enriched with all the quasi-simples mapping to each of its components is

$$\xymatrix{S[-1]\ar[d]&S\ar[d]\ar@{=}[r]&S\ar[dr]&&S\ar[dl]&S\ar[d]\ar@{=}[l]&S[1]\ar[d]\\
X[-1]\ar[r]&Y\ar[rr]&&M\ar[rr]&&X\ar[r]^{f}&Y[1]\\
T[-1]\ar[u]&T'\ar[u]\ar@{=}[r]&T'\ar[ur]&&T\ar[ul]&T\ar[u]\ar@{=}[l]&T'[1]\ar[u]
}$$

Since $X \in \St$ and $X$ is not isomorphic to either $Q_{S}^1$ or $Q_{S}^2$, Lemma \ref{lemma:StPoset} implies neither $Q_{S}^1$ nor $Q_{S}^2$ can map to $X$, and thus cannot be direct summands of $M$ (using Lemma \ref{lemma:RingelTriangles} to see that this would indeed induce a non-zero map from $Q_{S}^1$ or $Q_{S}^2$ to $X$) it follows from Proposition \ref{proposition:ProbingObjects}(\ref{proposition:ProbingObjects:Kardinality}) that $M$ is the direct sum of exactly two indecomposable objects, $M_{1}$ and $M_{2}$.  

It is now easy to see that we may assume $\phic(M_{1}) = \{S,T\}$ and $\phic(M_{1}) = \{S,T' \}$.  Proposition \ref{proposition:ProbingObjects}(\ref{proposition:ProbingObjects:Fiber}) now yields $M_{1} \cong X$ and $M_{2} \cong Y$, hence $M=X\oplus Y$, and thus $f=0$.

We now consider the latter case where neither $X$ nor $Y$ are the peripheral objects $Q_{S}^1$ or $Q_{S}^2$ and the composition $\xymatrix@1{T\ar[r]&X\ar[r]^{f}&Y[1]}$ is non-zero.  This yields
$$\xymatrix{S[-1]\ar[d]&S\ar[d]\ar@{=}[r]&S\ar[dr]&&S\ar[dl]&S\ar@{=}[l]\ar[d]&S[1]\ar[d]\\
X[-1]\ar[r]&Y\ar[rr]&&M\ar[rr]&&X\ar[r]^{f}&Y[1]\\
T[-1]\ar[u]\ar@{=}[r]&T'\ar[u]&&&&T\ar[u]\ar@{=}[r]&T'[1]\ar[u]
}$$
which is easily seen to be false since this would imply that either $Q_{S}^1$ or $Q_{S}^2$ would map to $X$.  Indeed, these are the only elements that do not have two different quasi-simples mapping to them and as such the only possible direct summands of $M$.
\end{proof}

\begin{lemma}\label{lemma:KeyLemma}
If $S$ is a quasi-simple object of $\CC$, then $\CC \cong \Db(\rep(\St^{\circ}))$.
\end{lemma}

\begin{proof}
First, due to Lemma \ref{lemma:StPTS} and Theorem \ref{theorem:PTS}, we may consider a full and exact embedding $\Db(\rep(\St^{\circ})) \to \CC$.  To show this is an equivalence of triangulated categories, we need to check that it is essentially surjective.  We will proceed in three steps.
\begin{enumerate}
\item
We start by showing that every indecomposable $Y$ with $\Hom(S,Y) \not= 0$ lies within the subcategory $\Db(\rep(\St^{\circ}))$ of $\CC$.  If $S$ does not map to a $\ZDi$-component, this follows directly from the definition of $\St$.  If $S$ does map to peripheral objects, $Q_{S}^1$ and $Q_{S}^2$, of a $\ZDi$-component, then we show there are non-zero morphisms $f_1:Q_{S}^1\to Y$ and $f_2:Q^2_{S}\to Y$.  First note that there is no map $Y \to Q_{S}^1$ since otherwise $Y \in \St$.  Now, the existence of $f$ follows from the diagram
$$\xymatrix{Q^1_{S}\ar[rd]^{f_1}\\
S\ar[u]\ar[r]&Y\ar[r]^{h}&C\ar[r]&S[1]}$$
where the bottom line is a triangle, $h:Y \to C$, is irreducible and the morphism $f_1:Q_{S}^1 \to Y$ is given by Lemma \ref{lemma:irredFact}.  Analogously, one proves the existence of $f_2:Q_{S}^2\to Y$.

Consider the following triangle, enriched with the quasi-simples mapping to each of its entries
$$\xymatrix{S[-1]\ar[d]&S\ar[d]&S\ar[dr]&&S\ar[dl]&S\ar[d]&S[1]\ar[d]\\
Y[-1]\ar[r]&X\ar[rr]_-{\begin{pmatrix}g_1\\g_2\end{pmatrix}}&&Q^1_{S}\oplus Q^2_{S} \ar[rr]_-{\begin{pmatrix}f_1&f_2\end{pmatrix}}&&Y\ar[r]&X[1]\\
T[-1]\ar[u]&T[-1]\ar[u]&&&&T\ar[u]&T\ar[u]
}$$
Since neither $f_1$ nor $f_2$ are zero, Lemma \ref{lemma:RingelTriangles} yields that $g$ and $g'$ are no isomorphisms.  In particular, $X$ cannot contain $Q_{S}^1$ or $Q_{S}^2$ as direct summands.  Due to the quasi-simple objects mapping to $X$, we may easily deduce that $X$ is an indecomposable object and $X \in \St$.  Since $\Db(\rep(\St^{\circ}))$ is an exact subcategory of $\CC$, we may conclude that $Y \in \Db(\rep(\St^{\circ}))$.
\item
We will now consider the more general case where $Y$ is an indecomposable object of $\CC$ such that $\Hom(S,Y[z])=0$ for all $z \in \bZ$.  Since the category $\CC$ is connected we may assume, without loss of generality, the existence of at least one indecomposable object $X$ of $\Db(\rep(\St^{\circ}))$ such that $\Hom(X,Y)\not=0$ or $\Hom(Y,X)\not=0$.  First, assume the former.  Since $\Db(\rep(\St^{\circ}))$ is generated by elements of $\St$ by taking finitely many cones and shifts, we may assume the existence of an indecomposable object $P \in \St$ such that $\Hom(P,Y)\not=0$.  Consider the triangle
$$\tri{E}{P\otimes\Hom(P,Y)}{Y}$$
as in Lemma \ref{lemma:CanonicalMorphism} where it has been proved that $\add E \cup \{P\}$ is a partial tilting set and $Y$ lies within the subcategory of $\CC$ generated by $E \oplus P$.  Let $E_{1}$ be any indecomposable direct summand of $E$ and consider the following triangle
$$\tri{E_{1}}{P}{C_{1}}.$$
Due to Lemma \ref{lemma:indecomposableCone} we may conclude that $C_{1}$ is indecomposable.  Applying the functor $\Hom(S,-)$ shows that either $\Hom(S,E_{1})\not=0$ or $\Hom(S,C_{1})\not=0$, and as such, either $E_{1} \in \Db(\rep(\St^{\circ}))$ or $C_{1} \in \Db(\rep(\St^{\circ}))$.  In both cases, since $\Db(\rep(\St^{\circ}))$ is an exact subcategory of $\CC$, we may conclude that $E_{1} \in \Db(\rep(\St^{\circ}))$, and thus that $E$ and hence also $Y$ lie within $\Db(\rep(\St^{\circ}))$.
\item
Finally, if $\Hom(Y,X)\not=0$ where $X$ is an indecomposable of $\Db(\rep(\St^{\circ}))$, consider the triangle $\tri{Y}{X}{C}$.  Due to Lemma \ref{lemma:RingelTriangles} we know there to be non-zero morphisms from $X$ to every direct summand of $C$.  In this case, it has been established in the second part of this proof that every indecomposable summand of $C$ lies in $\Ob \Db(\rep(\St^{\circ}))$ and hence also that $C \in \Ob \Db(\rep(\St^{\circ}))$.  Due to the fact that $\Db(\rep(\St^{\circ}))$ is an exact subcategory of $\CC$, we may conclude $Y \in \Db(\rep(\St^{\circ}))$.  This proves the assertion.
\end{enumerate}
\end{proof}

In previous lemma, we have used this easy lemma.

\begin{lemma}\label{lemma:indecomposableCone}
Let $X,Y \in \Ob \CC$ be indecomposables objects, and let $\tri{X}{Y}{Z}$ be the triangle built on a non-zero morphism $X\to Y$.  If $\Hom(Y,X[1])=0$ then $Z$ is indecomposable.
\end{lemma}

\begin{proof}
Using Lemma \ref{lemma:RingelTriangles} it should be clear that $Z$ has at most $\dim\Hom(Y,Z)$ direct summands.  Applying the functor $\Hom(Y,-)$ to the triangle $\tri{X}{Y}{Z}$, and using that $\dim\Hom(Y,Y)=1$ and $\Hom(Y,X[1])=0$, it follows easily that $\dim\Hom(Y,Z)=1$ and thus that $Z$ is indecomposable.
\end{proof}

We now prove our main result.

\begin{theorem}\label{theorem:MainTheorem}
A connected directed hereditary abelian $k$-linear Ext-finite category $\AA$ satisfying Serre duality is derived equivalent to $\repc(\PP)$ where $\PP$ is either a Dynkin quiver, $\AL$, or $\DL$ where $\LL$ is a locally discrete linearly ordered set without maximum or minimum.
\end{theorem}

\begin{proof}
From Theorem \ref{theorem:Components} we know that the only components of the Auslander-Reiten quiver of $\Db\AA$ are of the form $\ZAi$, $\ZAii$, $\ZDi$ or $\bZ Q$, where $Q$ is a quiver of Dynkin type.

First, assume $\CC$ has an Auslander-Reiten component of the form $\bZ Q$ where $Q$ is a Dynkin quiver.  It has been proven in Proposition \ref{proposition:ZDynkin} that $\Db\AA \cong \Db(\repc(Q^\circ))$

Thus we may now assume that the Auslander-Reiten quiver of $\Db\AA$ does not have a component of the form $\bZ Q$ where $Q$ is a Dynkin quiver.  The only possible components then are of the form $\ZAi$, $\ZAii$ or $\ZDi$.  It follows from Proposition \ref{proposition:ProbingComponents} that there is at least one wing $\WW$.  Fix a quasi-simple $S \in \WW$ and consider the set $\St$.  Lemma \ref{lemma:KeyLemma} yields that $\Db\AA \cong \Db(\rep(\St^{\circ}))$.

From Lemma \ref{lemma:StPoset} we know that $\St \cong A_{\LL'}^{\circ}$ or $\St \cong D_{\LL'}^{\circ}$ where $\LL' = (\bN)\cdot (\TT'\stackrel{\to}{\times}\bZ) \cdot (-\bN)$ for a certain linearly ordered set $\TT'$.  In the first case, we will give a poset $\AL$ such that $\Db \AA$ is equivalent to $\Db(\repc \AL)$; in the second case we will give a poset $\DL$ such that $\Db \AA$ is equivalent to $\Db(\repc \DL)$.

Thus, first assume $\St \cong A_{\LL'}^{\circ}$.  We will now consider the category $\Db(\rep A_{\LL})$ where $\LL = \TXZ$ and $\TT = \TT' \cdot \{ \ast \}$.  This category has already been discussed in \S\ref{section:AL}.  Fix the quasi-simple $T = A_{(*,0),(*,0)}$ of $\Db(\rep A_{\LL})$.  We may characterize $\Tt$ as 
$$\Tt = \{A_{(*,-n),(*,0)} \mid n \in \bN \} \cdot
\{A_{(t,z),(*,0)} \mid (t,z) \in \TT' \stackrel{\rightarrow}{\times} \bZ \} \cdot
\{A_{(*,1),(*,n+1)}[1] \mid n \in \bN \}$$
and thus $\Tt \cong A_{\LL'}^{\circ}$ as drawn in the following figure.
\begin{center}
		\psfrag{S}[][]{$S$}
		\psfrag{E}[][]{$\tau S[1]$}
		\psfrag{A}[][]{$\AA$}
		\psfrag{a1}[][]{$\AA[1]$}
		\includegraphics{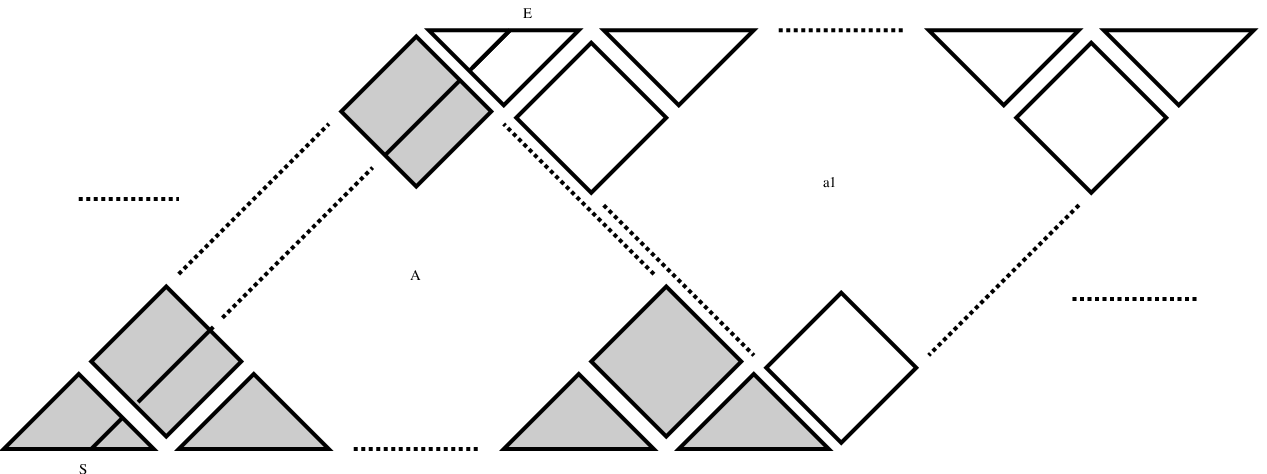}
\end{center}
By Lemma \ref{lemma:KeyLemma} we have that $\Db(\rep A_{\LL}) \cong \Db(\rep(A_{\LL'}))$, thus $\Db\AA \cong \Db(\repc \AL)$.

We will now consider the second case where $\St \cong D_{\LL'}^{\circ}$.  Consider $\Db(\repc(\DL))$ where $\LL = \TXZ$ for $\TT = \TT' \cdot \{ \ast \}$.  This category has already been discussed in \S\ref{section:DL}.  Fix the quasi-simple $T = A_{(*,0),(*,0)}$ of $\Db(\repc(\DL))$.  We may write $\Tt$ as 
$$\Tt = \{A_{(*,-n),(*,0)} \mid n \in \bN \} \cdot
\{A_{(t,z),(*,0)}  \mid (t,z) \in \TT' \stackrel{\rightarrow}{\times} \bZ \} \cdot
\{B_{(*,0),(*,n+1)} \mid n \in \bN \} \cdot
\{Q_{T}^1,Q_{T}^2\}$$
and thus
$\Tt \cong D_{\LL'}^{\circ}$ with $\LL' = \TT' \stackrel{\rightarrow}{\times} \bZ$ as in the following figure.
\begin{center}
		\psfrag{S}[][]{$S$}
		\psfrag{A}[][]{$\AA$}
		\psfrag{A1}[][]{$\AA[1]$}
		\includegraphics{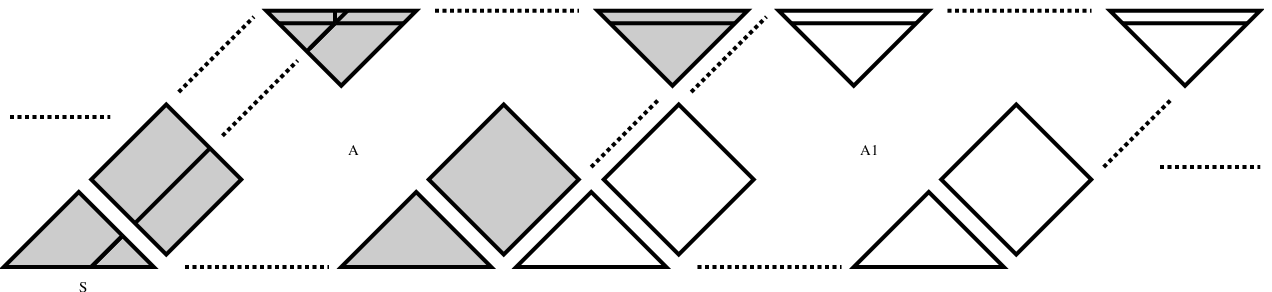}
\end{center}
  Lemma \ref{lemma:KeyLemma} yields $\Db(\repc(\DL)) \cong \Db(\rep(D_{\LL'}))$, and thus $\Db\AA \cong \Db(\repc(\DL))$.
\end{proof}

\begin{remark}
It follows from the proof of Theorem \ref{theorem:MainTheorem} that $\LL$ can also be chosen to be a \emph{bounded} linearly ordered set (i.e. $\LL$ has a maximal and a minimal element).
It then follows in particular that $\AA$ is derived equivalent to a hereditary category
which has both enough projectives and injectives.
\end{remark}

\begin{remark}
It is proven in Proposition \ref{proposition:DEquivalentToDynkin} that if, in the statement of Theorem \ref{theorem:MainTheorem}, $\PP$ is a Dynkin quiver, then $\AA$ is equivalent to $\mod kQ$ for a certain Dynkin quiver $Q$.
\end{remark}

\bibliographystyle{amsplain}

\providecommand{\bysame}{\leavevmode\hbox to3em{\hrulefill}\thinspace}
\providecommand{\MR}{\relax\ifhmode\unskip\space\fi MR }
\providecommand{\MRhref}[2]{%
  \href{http://www.ams.org/mathscinet-getitem?mr=#1}{#2}
}
\providecommand{\href}[2]{#2}

\end{document}